\documentclass[11pt]{article}
\usepackage{amsfonts,amsmath,amssymb,accents}
\usepackage{amsfonts,amsmath}
\usepackage{hyperref}
\usepackage{graphicx}
\usepackage{tikz}
\usetikzlibrary{calc,arrows,positioning,fit,petri}

\newcommand{\status}{}

\newcommand{\file}{}

\newcommand{\detail}[1]{\par\noi{\bf [Proof detail\ }{#1}
\hfill{\bf ]}\par\noi\hspace{-4pt}}
\renewcommand{\detail}[1]{}

\newcommand{\dis}{\displaystyle}

\newcommand{\tbf}{\textbf}
\newcommand{\pbf}[1]{{\rm\textbf{#1}}}

\newcommand{\med}{\medskip}
\newcommand{\noi}{\noindent}

\newcommand{\halmos}{\rule{1ex}{1.4ex}}
\def \qed {\nopagebreak{\hspace*{\fill}$\halmos$\medskip}}
\newcommand{\quand}{\quad\mbox{and}\quad}

\newtheorem{theorem}{Theorem}
\newtheorem{proposition}{Proposition}[section]
\newtheorem{corollary}[proposition]{Corollary}
\newtheorem{atheorem}[proposition]{Theorem}
\newtheorem{conjecture}[proposition]{Conjecture}
\newtheorem{lemma}[proposition]{Lemma}

\newtheorem{remark}[proposition]{Remark}

\newcommand{\bt}{\begin{theorem}}
\newcommand{\et}{\end{theorem}}
\newcommand{\bl}{\begin{lemma}}
\newcommand{\el}{\end{lemma}}
\newcommand{\bp}{\begin{proposition}}
\newcommand{\ep}{\end{proposition}}
\newcommand{\bcor}{\begin{corollary}}
\newcommand{\ecor}{\end{corollary}}
\newcommand{\br}{\begin{remark}\rm}
\newcommand{\er}{\end{remark}}
\newcommand{\bcon}{\begin{conjecture}}
\newcommand{\econ}{\end{conjecture}}

\newcommand{\be}{\begin{equation}}
\newcommand{\ee}{\end{equation}}
\newcommand{\ba}{\begin{array}}
\newcommand{\ea}{\end{array}}
\newcommand{\bc}{\be\begin{array}{r@{\,}c@{\,}l}}
\newcommand{\ec}{\end{array}\ee}

\newcommand{\al}{\alpha}

\newcommand{\ga}{\gamma}
\newcommand{\Ga}{\Gamma}
\newcommand{\de}{\delta}
\newcommand{\De}{\Delta}
\newcommand{\eps}{\varepsilon}
\newcommand{\la}{\lambda}
\newcommand{\La}{\Lambda}
\newcommand{\sig}{\sigma}

\newcommand{\tet}{\theta}
\newcommand{\Om}{\Omega}
\newcommand{\om}{\omega}

\newcommand{\si}{\ensuremath{\sigma}}

\newcommand{\Ai}{{\cal A}}

\newcommand{\Ci}{{\cal C}}
\newcommand{\Di}{{\cal D}}
\newcommand{\Ei}{{\cal E}}
\newcommand{\Fi}{{\cal F}}

\newcommand{\Pc}{{\cal P}}

\newcommand{\Si}{{\cal S}}

\newcommand{\R}{{\mathbb R}}
\newcommand{\N}{{\mathbb N}}
\newcommand{\Z}{{\mathbb Z}}

\newcommand{\E}{{\mathbb E}}
\renewcommand{\P}{{\mathbb P}}

\newcommand{\volgt}{\ensuremath{\Rightarrow}}
\newcommand{\up}{\uparrow}
\newcommand{\down}{\downarrow}
\newcommand{\sub}{\subset}
\newcommand{\beh}{\backslash}

\newcommand{\asto}[1]{\underset{{#1}\to\infty}{\longrightarrow}}
\newcommand{\Asto}[1]{\underset{{#1}\to\infty}{\Longrightarrow}}

\newcommand{\ti}{\tilde}
\newcommand{\dgg}{\dagger}
\newcommand{\ov}{\overline}

\newcommand{\subb}[2]{_{\ba{c}\scriptstyle{#1}\\[-.15cm]\scriptstyle{#2}\ea}}

\newcommand{\ffrac}[2]{{\textstyle\frac{{#1}}{{#2}}}}
\newcommand{\dif}[1]{\ffrac{\partial}{\partial{#1}}}

\newcommand{\di}{\mathrm{d}}
\newcommand{\half}{{[0,\infty)}}
\newcommand{\expo}{\mbox{\large\it e}}
\newcommand{\ex}[1]{\expo^{\,\textstyle{#1}}}

\newcommand{\cov}{{\rm Cov}}

\newcommand{\poi}{\om}

\newcommand{\Pfp}{\Pc_{{\rm fin},\,+}}
\newcommand{\Pf}{\Pc_{\rm fin}}
\newcommand{\tiPfp}{\ti\Pc_{{\rm fin},\,+}}
\newcommand{\acPfp}{\Pc'_{{\rm fin},\,+}}
\newcommand{\tiPf}{\ti\Pc_{\rm fin}}

\newcommand{\hook}[1]{\accentset{{#1}}{\hookrightarrow}}
\newcommand{\nucirc}{{\accentset{\circ}{\nu}}}
\newcommand{\Pfi}[1]{\Pc_{{\rm fin},\,{#1}}}
\newcommand{\captimes}{\tikz[baseline=-3.3pt]{\node at (0,0) {$\cap$};
\node at (0,0) {$\scriptscriptstyle\times$};}}
\newcommand{\lla}{\langle\!\langle}
\newcommand{\rra}{\rangle\!\rangle}

\begin{document}

\makeatletter\@addtoreset{equation}{section}
\makeatother\def\theequation{\thesection.\arabic{equation}} 

\renewcommand{\labelenumi}{{\rm (\roman{enumi})}}

\title{\vspace*{-3cm}Subcritical contact processes seen from a
typical infected site}
\author{}
\author{Anja Sturm
\and
Jan~M.~Swart}

\date{{\small\file}\today}

\maketitle\vspace{-.8cm}
\status

\begin{abstract}\noi
What is the long-time behavior of the law of a contact process started with a
single infected site, distributed according to counting measure on the
lattice? This question is related to the configuration as seen from a typical
infected site and gives rise to the definition of so-called eigenmeasures,
which are possibly infinite measures on the set of nonempty configurations
that are preserved under the dynamics up to a multiplicative constant. In this
paper, we study eigenmeasures of contact processes on general countable groups
in the subcritical regime. We prove that in this regime, the process has a
unique spatially homogeneous eigenmeasure. As an application, we show that the
exponential growth rate is continuously differentiable and strictly decreasing
as a function of the recovery rate, and we give a formula for the derivative in
terms of the eigenmeasures of the contact process and its dual.
\end{abstract}
\vspace{.5cm}

\noi
{\it Acknowledgement.} Work sponsored by GA\v CR grant 201/09/1931.

%
%


{\setlength{\parskip}{-2pt}\tableofcontents}

\newpage

\section{Introduction and main results}

\subsection{Introduction}\label{S:intro}

It is known that contact processes on regular trees behave quite differently
from contact processes on the $d$-dimensional integer lattice $\Z^d$. Indeed,
if $\la_{\rm c}$ and $\la'_{\rm c}$ denote the critical infection rates
associated with global and local survival, respectively, then one has
$\la_{\rm c}<\la'_{\rm c}$ on trees while $\la_{\rm c}=\la'_{\rm c}$ on
$\Z^d$. For $\la>\la'_{\rm c}$, the process exhibits complete convergence and
the upper invariant law is the only nontrivial invariant law, while on trees,
in the intermediate regime $\la_{\rm c}<\la\leq\la'_{\rm c}$, there is a
multitude of (not spatially homogeneous) invariant laws. The situation is
reminiscent of what is known about unoriented percolation on transitive
graphs, where one has uniqueness of the infinite cluster if the graph is
amenable, while it is conjectured, and proved in some cases, that on
nonamenable graphs there is an intermediate parameter regime with infinitely
many infinite clusters. We refer to \cite{Lig99} as a general reference to
contact processes on $\Z^d$ and trees and \cite{Hag11} for percolation beyond
$\Z^d$.

In general, it is not hard (but also not very interesting) to determine the
limit behavior of contact processes started from a spatially homogeneous
(i.e., translation invariant) initial law. On the other hand, it seems much
more difficult to study the process started with a finite number of infected
sites. For example, the proof of complete convergence on $\Z^d$
 \cite[Section~I.2]{Lig99} is rather involved and cannot easily be adapted to
 other lattices. In a similar vein, we mention that it is an open problem to
 prove that $\la_{\rm c}=\la'_{\rm c}$ on any amenable transitive graph; this
 seems to be quite hard. As an intermediate problem, in \cite[Problem~1 from
 Section~1.5]{Swa09}, it has been proposed to study the process started with
a single infected site, chosen uniformly from the lattice. For infinite
lattices, the resulting `law' at time $t$ will be an infinite
measure. However, as shown in \cite[Lemma~4.2]{Swa09}, conditioning such a
measure on the origin being infected yields a probability law, which can be
interpreted as the process seen from a typical infected site.

There is a close connection between the law of the process seen from a typical
infected site and the exponential growth rate $r$ of 
the expected number of infected sites of a contact process. This can be
understood by realizing that the number of healthy sites surrounding a typical
infected site determines the number of infections that can be made and hence
the speed at which the infection grows. In the context of infinite laws, which
cannot be normalized, it is natural to generalize the concept of an invariant
measure to an `eigenmeasure', which is a measure on the set of nonempty
configurations that is preserved under time evolution up to a multiplicative
constant. Alternatively, such eigenmeasures can be thought of as the
equivalent of a quasi-stationary law (as introduced in \cite{DS67}) in the
setting of interacting particle systems. In particular, if the suitably
rescaled `law' at time $t$ of the process started with a single, uniformly
distributed site has a nontrivial long-time limit, then it follows from
results in \cite{Swa09} that such a limit `law' must be an eigenmeasure whose
eigenvalue is the exponential growth rate $r$ of the process.

In the present paper, we study eigenmeasures of subcritical contact processes
on general countable groups. Our set-up includes translation-invariant contact
processes on $\Z^d$ and on regular trees, as well as long-range processes and
asymmetric processes. We will show that such processes have a unique
homogeneous eigenmeasure which is the vague limit of the rescaled law at time
$t$ of the process started in any homogeneous, possibly infinite, initial
law. As an application of our results, we give an expression for the
derivative of the exponential growth rate as a function of the recovery rate
in terms of the eigenmeasures of the process and its dual, and we use this to
show that this derivative is strictly negative and continuous.

\subsection{Contact processes on groups}\label{S:group}

We need to define the class of contact processes that we will be interested
in, fix notation, and recall some well-known facts.
Let $\La$ be a finite or countably infinite group with group action
$(i,j)\mapsto ij$, inverse operation $i\mapsto i^{-1}$, and unit element $0$
(also refered to as the origin). Let $a:\La\times\La\to\half$ be a function
such that $a(i,i)=0$ $(i\in\La)$ and
\be\ba{rl}\label{assum}
{\rm (i)}&a(i,j)=a(ki,kj)\qquad\qquad(i,j,k\in\La),\\[5pt]
{\rm (ii)}&\dis|a|:=\sum_{i\in \La}a(0,i)<\infty,\\[5pt]
\ec
and let $\de\geq 0$. By definition, the \emph{$(\La,a,\de)$-contact process} is
the Markov process $\eta=(\eta_t)_{t\geq 0}$, taking values in the space
$\Pc=\Pc(\La):=\{A:A\sub\La\}$ consisting of all subsets of $\La$, with the
formal generator
\bc\label{Gdef}
Gf(A)&:=&\dis\sum_{i,j\in\La}a(i,j)1_{\{i\in A\}}1_{\{j\notin A\}}
\{f(A\cup\{j\})-f(A)\}\\[5pt]
&&\dis+\de\sum_{i\in\La}1_{\{i\in A\}}\{f(A\beh\{i\})-f(A)\}.
\ec
If $i\in\eta_t$, then we say that the site $i$ is infected at time $t$;
otherwise it is healthy. Then (\ref{Gdef}) says that an infected site $i$
infects another site $j$ with \emph{infection rate} $a(i,j)\geq 0$, and infected
sites become healthy with \emph{recovery rate} $\de\geq 0$.

We will usually assume that the infection rates are irreducible in some sense
or another. To make this precise, let us write $i\hook{a} j$ if the site $j$
can be infected through a chain of infections starting from $i$. Then we say
that $a$ is \emph{irreducible} if $i\hook{a} j$ for all $i,j\in\La$.
Equivalently, this says that for all $\La'\sub\La$ with
$\La'\neq\emptyset,\La$, there exist $i\in\La'$ and $j\in\La\beh\La'$ such
that $a(i,j)>0$. Similarly, we say that $a$ is \emph{weakly irreducible} if for
all $\La'\sub\La$ with $\La'\neq\emptyset,\La$, there exist $i\in\La'$ and
$j\in\La\beh\La'$ such that $a(i,j)\vee a(j,i)>0$.
Finally, we will sometimes need the intermediate condition
\be\label{irr}
\forall i,j\in\La:\ \exists k,l\in\La:
\  k\hook{a} i,\  k\hook{a} j,\ i\hook{a} l,\ j\hook{a} l.
\ee
In words, this says that for any two sites $i,j$ there exists a site $k$ from
which both $i$ and $j$ can be infected, and a site $l$ that can be infected
both from $i$ and from $j$.  If the rates $a$ are symmetric, or more generally
if one has $a(i,j)>0$ iff $a(j,i)>0$, then all three conditions are
equivalent. In general, irreducibility implies (\ref{irr}) which implies weak
irreducibility, but none of the converse implications hold.

It is well-known that contact processes can be constructed by a graphical
representation. Let $\poi=(\poi^{\rm r},\poi^{\rm i})$ be a pair of
independent, locally finite random subsets of $\La\times\R$ and
$\La\times\La\times\R$, respectively, produced by Poisson point processes with
intensity $\de$ and $a(i,j)$, respectively. This is usually visualized by
plotting $\La$ horizontally and $\R$ vertically, marking points
$(i,s)\in\poi^{\rm r}$ with a recovery symbol (e.g., $\ast$), and drawing an
infection arrow from $(i,t)$ to $(j,t)$ for each $(i,j,t)\in\poi^{\rm i}$. For
any $(i,s),(j,u)\in\La\times\R$ with $s\leq u$, by definition, an \emph{open
path} from $(i,s)$ to $(j,u)$ is a cadlag function $\pi:[s,u]\to\La$ such
that $\{(\pi(t),t):t\in[s,u]\}\cap\poi^{\rm r}=\emptyset$ and
$(\pi(t-),\pi(t),t)\in\poi^{\rm i}$ whenever $\pi(t-)\neq\pi(t)$. Thus, open
paths must avoid recovery symbols and may follow infection arrows. We write
$(i,s)\leadsto(j,u)$ to indicate the presence of an open path from $(i,s)$ to
$(j,u)$. Then, for any $s\in\R$, we can construct a $(\La,a,\de)$-contact
process started in an initial state $A\in\Pc$ by setting
\be\label{etaA}
\eta^{A,s}_t:=\{j\in\La:(i,s)\leadsto(j,s+t)\mbox{ for some }i\in A\}
\qquad(A\in\Pc,\ s\in\R,\ t\geq 0).
\ee
In particular, we set $\eta^A_t:=\eta^{A,0}_t$. Note that this construction
defines contact processes with different initial states on the same
probability space, i.e., the graphical representation provides a natural
coupling between such processes. Moreover, the graphical representation shows
that the contact process is essentially a sort of oriented percolation model
(in continuous time but discrete space).

Since the graphical representation is also defined for negative times we can, in
analogy to (\ref{etaA}), define `backward' or `dual' processes by
\be
\eta^{\dgg\,A,s}_t:=\{j\in\La:(j,s-t)\leadsto(i,s)\mbox{ for some }i\in A\}
\qquad(A\in\Pc,\ s\in\R,\ t\geq 0).
\ee
In particular, we set $\eta^{\dgg\,A}_t:=\eta^{\dgg\,A,0}_t$. It is not hard
to see that $(\eta^{\dgg\,A,s}_t)_{t\geq 0}$ is a $(\La,a^\dgg,\de)$-contact
process, where we define \emph{reversed infection rates} as
$a^\dgg(i,j):=a(j,i)$. Since
\be
\big\{\eta^A_t\cap B\neq\emptyset\big\}
=\big\{(i,0)\leadsto(j,t)\mbox{ for some }i\in A,\ j\in B\big\}
=\big\{\eta^A_0\cap\eta^{\dgg\,B,t}_t\neq\emptyset\big\}
\qquad(0\leq s\leq t)
\ee
and the process $\eta^{\dgg\,B,t}$ is equal in law with $\eta^{\dgg\,B}$, we
see that the $(\La,a,\de)$-contact process and $(\La,a^\dgg,\de)$-contact
process are dual in the sense that
\be\label{dual}
\P[\eta^A_t\cap B\neq\emptyset]=\P[A\cap\eta^{\dgg\,B}_t\neq\emptyset]
\qquad(A,B\in\Pc,\ t\geq 0).
\ee
We note that unless $a=a^\dgg$ or the group $\La$ is abelian, the
$(\La,a,\de)$- and $(\La,a^\dgg,\de)$-contact processes have in general
different dynamics and need to be distinguished. (If $\La$ is abelian, then
the $(\La,a,\de)$- and $(\La,a^\dgg,\de)$-contact processes can be mapped into
each other by the transformation $i\mapsto i^{-1}$.) 

We say that the $(\La,a,\de)$-contact process \emph{survives} if
$\P[\eta^A_t\neq\emptyset\ \forall t\geq 0]>0$
for some, and hence for all nonempty $A$ of finite cardinality $|A|$.
We call
\be\label{dec}
\de_{\rm c}=\de_{\rm c}(\La,a):=\sup\big\{\de\geq 0:
\mbox{ the $(\La,a,\de)$-contact process survives}\big\}
\ee
the \emph{critical recovery rate}. It is known that $\de_{\rm c}<\infty$. If
$\La$ is finitely generated, then moreover $\de_{\rm c}>0$ provided $a$ is
weakly irreducible \cite[Lemma~4.18]{Swa07}, but for non-finitely generated
groups irreducibility is in general not enough to guarantee $\de_{\rm c}>0$
\cite{AS10}. It is well-known that
\be
\P\big[\eta_t^\La \in\cdot\,\big]\Asto{t}\ov\nu,
\ee
where $\ov\nu$ is an invariant law of the $(\La,a,\de)$-contact process, known
as the \emph{upper invariant law}. Using duality, it is not hard to prove that
$\ov\nu=\de_\emptyset$ if the $(\La,a^\dgg,\de)$-contact process dies out,
while $\ov\nu$ is concentrated on the nonempty subsets of $\La$ if the process
survives.

It follows from subadditivity (see \cite[Lemma~1.1]{Swa09}) that any
$(\La,a,\de)$-contact process has a well-defined exponential growth rate,
i.e., there exists a constant $r=r(\La,a,\de)$ with $-\de\leq r\leq|a|-\de$
such that
\be\label{rdef}
r=\lim_{t\to\infty}\ffrac{1}{t}\log\E\big[|\eta^A_t|\big]\qquad(0<|A|<\infty).
\ee
In this article, we are concerned with subcritical contact processes for which $r<0.$
The following theorem lists some properties of the function $r(\La,a,\de)$.

\addtocounter{theorem}{-1}

\bt\pbf{(Properties of the exponential growth rate)}\label{T:rprop}\ \\
For any $(\La,a,\de)$-contact process:
\begin{itemize}
\item[\rm(a)] $r(\La,a,\de)=r(\La,a^\dgg,\de)$.
\item[\rm(b)] The function $\de\to r(\La,a,\de)$ is nonincreasing and
Lipschitz continuous on $\half$, with Lipschitz constant 1.
\item[\rm(c)] If $r(\La,a,\de)>0$, then the $(\La,a,\de)$-contact process
survives.
\item[\rm(d)] $\{\de\geq 0:r(\La,a,\de)<0\}=(\de_{\rm c},\infty)$.
\end{itemize}
\et

The (easy) proofs of parts~(a)--(c) can be found in
\cite[Theorem~1.2]{Swa09}. The analogue of part~(d) for unoriented percolation
on $\Z^d$ was first proved by Menshikov \cite{Men86} and Aizenman and Barsky
\cite{AB87}. Using the approach of the latter paper, Bezuidenhout and Grimmett
\cite[formula (1.13)]{BG91} proved the statement in part~(d) for contact
processes on $\Z^d$. This has been generalized to processes on general
transitive graphs in \cite{AJ07}. As we point out in Appendix~\ref{A:decay},
their arguments are not restricted to graphs but apply in the generality we
need here. We note that it follows from parts~(a) and (d) that $\de_{\rm
c}(\La,a)=\de_{\rm c}(\La,a^\dgg)$. In general, it is not known if survival
of a $(\La,a,\de)$-contact process implies survival of the dual
$(\La,a^\dgg,\de)$-contact process but any counterexample would have to be at
$\de=\de_{\rm c}$, while by \cite[Corollary~1.3]{Swa09}, $\La$ would have to
be amenable. If $\La$ is a finitely generated group of subexponential growth
and the infection rates satisfy an exponential moment condition (for example,
if $\La=\Z^d$ and $a$ is nearest-neighbor), then $r\leq 0$
\cite[Thm~1.2~(e)]{Swa09}, but in general (e.g.\ on trees), it is possible
that $r>0$. Indeed, one of the main results of \cite{Swa09} says that if $\La$
is nonamenable, the $(\La,a,\de)$-contact process survives, and the infection
rates satisfy the irreducibility condition (\ref{irr}), then $r>0$
\cite[Thm.~1.2~(f)]{Swa09}.

\subsection{Locally finite starting measures}\label{S:locfin}

We will be interested in the contact process started in initial `laws' that
are infinite measures. To do this properly, we need a bit of theory. Recall
that $\Pc=\Pc(\La)$ denotes the space of all subsets of $\La$. 
We let $\Pc_+:=\{A:|A|>0\}$ and $\Pf:=\{A:|A|<\infty\}$ denote the subspaces
consisting of all nonempty, respectively finite subsets of $\La$, and write
$\Pfp:=\Pf\cap\Pc_+$. We observe that $\Pc\cong\{0,1\}^{\La}$ and equip it
with the product topology and Borel-\si-field. 
Note that since $\Pc$ is compact, $\Pc_+=\Pc\beh\{\emptyset\}$ is a locally
compact space. Recall that a measure on a locally compact space is {\em
locally finite} if it gives finite mass to compact sets, and that a sequence
of locally finite measures converges vaguely if the integrals of all compactly
supported, continuous functions converge. We cite the following simple facts
from \cite[Lemmas~3.1 and 3.2]{Swa09}.

\bl\pbf{(Locally finite measures)}\label{L:locfin}
Let $\mu$ be a measure on $\Pc_+$. Then the following statements are
equivalent:
\begin{enumerate}
\item $\mu$ is locally finite.
\item $\int\mu(\di A)1_{\{i\in A\}}<\infty$ for all $i\in\La$. 
\item $\int\mu(\di A)1_{\{A\cap B\neq\emptyset\}}<\infty$ for all $B\in\Pfp$.
\end{enumerate}
Moreover, if $\mu_n,\mu$ are locally finite measures on $\Pc_+$, then 
the $\mu_n$ converge vaguely to $\mu$ if and only if
\be
\int\mu_n(\di A)1_{\{A\cap B\neq\emptyset\}}
\asto{n}\int\mu(\di A)1_{\{A\cap B\neq\emptyset\}}
\qquad(B\in\Pfp).
\ee
\el

We will sometimes deal with locally finite measures on $\Pc_+$ that are
concentrated on $\Pf$. We will refer to such measures as `locally
finite measures on $\Pfp$' (even though `locally finite' refers to the
topology on $\Pc_+$). For such measures, we will sometimes need another,
stronger form of convergence than vague convergence. For each $i\in\La$, we
define
\be\label{Pidef}
\Pc_i:=\{A\in\Pc:i\in A\}
\quad\mbox{and}\quad
\Pfi{i}:=\Pf\cap\Pc_i.
\ee
Note that $\Pfi{i}$ is a countable set. We let $\mu|_{\Pfi{i}}$ denote the
restriction of a measure $\mu$ to $\Pfi{i}$. If $\mu_n,\mu$ are locally finite
measures on $\Pfp$, then we say that the $\mu_n$ converge to $\mu$ {\em
locally on $\Pfp$}, if for each $i\in\La$, the $\mu_n|_{\Pfi{i}}$ converge
weakly to $\mu|_{\Pfi{i}}$ with respect to the discrete topology on
$\Pfi{i}$. It can be shown that local convergence on $\Pfp$ implies vague
convergence (see Proposition~\ref{P:locon} below), but the converse is not
true. For example, if $\La=\Z$, then using Lemma~\ref{L:locfin} it is not hard
to see that we have the vague convergence
\be\label{example}
\sum_{i\in\Z}\de_{\{i,i+n\}}=:\mu_n\Asto{n}
\mu:=2\sum_{i\in\Z}\de_{\{i\}},
\ee
(where $\de_A$ denotes the delta-measure at a point $A\in\Pc_+$) but the
$\mu_n$ do not converge locally on $\Pfp$.

We now turn our attention to contact processes started in infinite initial
`laws'. For a given $(\La,a,\de)$-contact process, we define subprobability
kernels $P_t$ $(t\geq 0)$ on $\Pc_+$ by
\be\label{Ptdef}
P_t(A,\,\cdot\,):=\P\big[\eta^A_t\in\cdot\,\big]\big|_{\Pc_+}
\qquad(t\geq 0,\ A\in\Pc_+),
\ee
where $|_{\Pc_+}$ denotes restriction to $\Pc_+$, and we define $P^\dgg_t$
similarly for the dual $(\La,a^\dgg,\de)$-contact process. For any measure
$\mu$ on $\Pc_+$, we write
\be
\label{muPtdef}
\mu P_t:=\int\mu(\di A)P_t(A,\,\cdot\,)\qquad(t\geq 0),
\ee
which is the restriction to $\Pc_+$ of the `law' at time $t$ of the
$(\La,a,\de)$-contact process started in the initial (possibly infinite) `law'
$\mu$.

For $A\sub\La$ and $i\in\La$, we write $iA:=\{ij:j\in A\}$, and for any
$\Ai\sub\Pc$ we write $i\Ai:=\{iA:A\in\Ai\}$. We say that a measure $\mu$ on
$\Pc$ is (spatially) \emph{homogeneous} if it is invariant under the left
action of the group, i.e., if $\mu(\Ai)=\mu(i\Ai)$ for each $i\in\La$ and 
measurable $\Ai\sub\Pc$. If $\mu$ is a homogeneous, locally finite measure on
$\Pc_+$, then $\mu P_t$ is a homogeneous, locally finite measure on $\Pc_+$
for each $t\geq 0$ (see \cite[Lemma~3.3]{Swa09} or Lemma~\ref{L:infstart}
below).

For processes started in homogeneous, locally finite measures, we have a
useful sort of analogue of the duality formula (\ref{dual}). To formulate
this, we need two more definitions. For any measure $\mu$ on $\Pc_+$,
we define
\be\label{cmudef}
\lla\mu\rra:=\int\mu(\di A)|A|^{-1}1_{\{0\in A\}},
\ee
where $|A|^{-1}:=0$ if $A$ is infinite. Note that if each set $A\in\Pfp$
carries mass $\mu(\{A\})$, and this mass is distributed evenly among all
points in $A$, then $\lla\mu\rra$ is the mass received at the origin.

Next, for any measures $\mu,\nu$ on $\Pc_+$, we let
$\mu\captimes\nu$ denote the restriction to $\Pc_+$ of the image of the
product measure $\mu\otimes\nu$ under the map $(A,B)\mapsto A\cap B$. Note
that
\be\label{psidef}
\int\mu\captimes\nu\,(\di C)f(C):=\int\mu(\di A)\int\nu(\di B)f\big(A\cap B)
\ee
for any bounded measurable $f:\Pc\to\R$ satisfying $f(\emptyset)=0$. We call
$\mu\captimes\nu$ the \emph{intersection measure} of $\mu$ and $\nu$.
It is not hard to show (see Lemma~\ref{L:intsect} below) that 
$\mu\captimes\nu$ is locally finite if $\mu$ and $\nu$ are. Note that
if $\mu$ and $\nu$ are probability measures, then $\mu\captimes\nu$ is the law
of the intersection of two independent random sets with laws $\mu$ and $\nu$,
restricted to the event that this intersection is nonempty. In particular,
normalizing $\mu\captimes\nu$ yields the conditional law given this event.

With these definitions, we have the following lemma, whose proof  can
be found in Section~\ref{S:infstart}.

\bl\pbf{(Duality for infinite initial laws)}\label{L:infdual}
Let $\mu,\nu$ be homogeneous, locally finite measures on $\Pc_+$. Then
\be\label{infdual}
\lla\mu P_t\captimes\nu\rra=\lla\mu\captimes\nu P^\dgg_t\rra
\qquad(t\geq 0),
\ee
and $\mu P_t\captimes\nu$ is concentrated on $\Pfp$ if and only if
$\mu\captimes\nu P^\dgg_t$ is.
\el

\noi
\tbf{Remark} If $|\mu|:=\mu(\Pc_+)$ denotes the total mass of a finite measure
on $\Pc_+$, then the duality formula (\ref{dual}) is easily seen to imply that
$|\mu P_t\captimes\nu|=|\mu\captimes\nu P^\dgg_t|$ for any finite measures
$\mu,\nu$ on $\Pc_+$. One can think of (\ref{infdual}) as an analogue of this
for infinite (but homogeneous) measures.

\subsection{Eigenmeasures}\label{S:eigen}

Following \cite{Swa09}, we say that a measure $\mu$ on $\Pc_+$ is an {\em
eigenmeasure} of the $(\La,a,\de)$-contact process if $\mu$ is nonzero,
locally finite, and there exists a constant $\la\in\R$ such that
\be\label{eigen}
\mu P_t=e^{\la t}\mu\qquad(t\geq 0).
\ee
We call $\la$ the associated \emph{eigenvalue}.

It follows from \cite[Prop.~1.4]{Swa09} that each $(\La,a,\de)$-contact
process has a homogeneous eigenmeasure $\nucirc$ with eigenvalue
$r=r(\La,a,\de)$. In general, it is not known if $\nucirc$ is (up to a
multiplicative constant) unique. Under the irreducibility condition
(\ref{irr}), it has been shown in \cite[Thm.~1.5]{Swa09} that if the upper
invariant measure $\ov\nu$ of a $(\La,a,\de)$-contact process is concentrated
on $\Pc_+$ and $r(\La,a,\de)=0$, then $\nucirc$ is unique up to a
multiplicative constant and in fact $\nucirc=c\,\ov\nu$ for some $c>0$. The
main aim of the present paper is to investigate eigenmeasures in the
subcritical case $r<0$. Here is our first main result.

\bt\pbf{(Eigenmeasures in the subcritical case)}\label{T:eigcon}
Assume that the infection rates sa\-tisfy the irreducibility condition
(\ref{irr}) and that the exponential growth rate from (\ref{rdef}) satisfies
$r<0$. Then there exist, up to multiplicative constants, unique
homogeneous eigenmeasures $\nucirc$ and $\nucirc^\dgg$ of the $(\La,a,\de)$-
and $(\La,a^\dgg,\de)$-contact processes, respectively. These
eigenmeasures have eigenvalue $r$ and are concentrated on $\Pf$. If $\mu$ is any
nonzero, homogeneous, locally finite measure on $\Pc_+$, then
\be\label{eeigcon}
e^{-rt}\mu P_t\Asto{t}c\,\nucirc,
\ee
where $\Rightarrow$ denotes vague convergence of locally
finite measures on $\Pc_+$ and $c>0$ is a constant, given by
\be\label{cform}
c=\frac{\lla\mu\captimes\nucirc^\dgg\rra}{\lla\nucirc\captimes\nucirc^\dgg\rra}.
\ee
If $\mu$ is concentrated on $\Pfp$, then (\ref{eeigcon}) holds in the sense of
local convergence on $\Pfp$.
\et
The proof of Theorem \ref{eeigcon} will be completed in Section
\ref{S:exunique}.\med

\noi
\tbf{Remark} Since $\nucirc$ and $\nucirc^\dgg$ are infinite measures, their
normalizations are somewhat arbitrary. For definiteness, we will usually adopt
the convention that $\int\nucirc(\di A)1_{\{0\in A\}}=1=\int\nucirc^\dgg(\di
A)1_{\{0\in A\}}$. Theorem~\ref{T:eigcon} holds regardless of the choice of
normalization.\med

\subsection{The process seen from a typical infected site}

We next set out to explain the connection of eigenmeasures and the process as
seen from a typical infected site, and formulate our second main result, which
gives a formula for the derivative of the exponential growth rate.

Let $(\eta^{\{0\}}_t)_{t\geq 0}$ be a $(\La,a,\de)$-contact process, started
with a single infected site at the origin, where
$\eta^{\{0\}}_t=\eta^{\{0\}}_t(\om)$ is defined on some underlying probability
space $(\Om,\Fi,\P)$. Then, for each $t\geq 0$, we can define a new
probability law $\hat\P_t$ on a suitably enriched probability space $\hat\Om$
that also contains a $\La$-valued random variable $\iota$, by setting
\be\label{Campbell}
\hat\P_t\big[\om\in\Ai,\ \iota=i\big]
:=\frac{\P[\om\in\Ai,\ i\in\eta^{\{0\}}_t(\om)]}{\E[|\eta^{\{0\}}_t|]}
\qquad(\Ai\in\Fi,\ i\in\La).
\ee
The law $\hat\P_t$ is a Campbell law (closely related to the more well-known
Palm laws). In words, $\hat\P_t$ is obtained from the original law $\P$ by
size-biasing on the number $|\eta^{\{0\}}_t|$ of infected sites at time $t$
and then choosing one site $\iota$ from $\eta^{\{0\}}_t$ with equal
probabilities.

Let $\mu_t:=\sum_{i\in\La}\P[\eta^{\{i\}}_t\in\cdot\,]|_{\Pc_+}$ be the
infinite `law' of the process started with a single infection at a uniformly
chosen site in the lattice. Then, defining conditional probabilities for
infinite measures in the natural way, it has been shown in
\cite[Lemma~4.2]{Swa09} that
\be\label{fromtypic}
\mu_t\big(\,\cdot\,\big|\,\{A:0\in A\}\big)
:= \frac{\mu_t\big(\,\cdot\,\cap\,\{A:0\in A\}\big)}
{\mu_t\big( \{A:0\in A\}\big)}
=\hat\P_t\big[\iota^{-1}\eta^{\{0\}}_t\in\cdot\,\big],
\ee
i.e., $\mu_t$ conditioned on the origin being infected describes the
distribution of $\eta^{\{0\}}_t$ under the Campbell law $\hat\P_t$ with the
`typical infected site' $\iota$ shifted to the origin. 

In view of this, Theorem~\ref{T:eigcon} gives information about the long-time
limit law of the process seen from a typical infected site. Indeed, it is easy
to see that Theorem~\ref{T:eigcon} implies the weak convergence of the
probability measures in (\ref{fromtypic}) to $\nucirc(\,\cdot\,\big|\,\{A:0\in
A\})$.


To see the connection of this with the derivative of the exponential growth
rate, let $\eta^{\de,\,\{0\}}_t$ denote the
process with a given recovery rate $\de$ (and $(\La,a)$ fixed), constructed
with the graphical representation. A version of Russo's formula (see
\cite[formula~(3.10)]{Swa09} and compare \cite[Thm~2.25]{Gri99}) tells us that
\be\label{Russo}
-\dif{\de}\frac{1}{t}\log\E\big[|\eta^{\de,\,\{0\}}_t|\big]
=\frac{1}{t}\int_0^t
\hat\P_t\big[\exists j\in\La\mbox{ s.t.\ }
(0,0)\leadsto_{(j,s)}(\iota,t)\big]\di s,
\ee
where $(0,0)\leadsto_{(j,s)}(\iota,t)$ denotes the event that in the graphical
representation, all open paths from $(0,0)$ to $(\iota,t)$ lead through
$(j,s)$. In other words, the right-hand side of (\ref{Russo}) is the fraction
of time that there is a \emph{pivotal} site on the way from $(0,0)$ to the
typical site $(\iota,t)$.

By grace of Theorem~\ref{T:eigcon}, we are able to control the long-time
limit of formula (\ref{Russo}), leading to the following result, whose proof
will be completed at the end of Section \ref{S:deriv}.

\bt\pbf{(Derivative of the exponential growth rate)}\label{T:difr}
Assume that the infection rates satisfy the irreducibility condition
(\ref{irr}). For $\de\in(\de_{\rm c},\infty)$, let $\nucirc_\de$ and
$\nucirc^\dgg_\de$ denote the homogeneous eigenmeasures of the $(\La,a,\de)$-
and $(\La,a^\dgg,\de)$-contact processes, respectively, normalized such that
$\int\nucirc_\de(\di A)1_{\{0\in A\}}=1=\int\nucirc^\dgg_\de(\di A)1_{\{0\in
A\}}$. Then the map $(\de_{\rm c},\infty)\ni\de\mapsto\nucirc_\de$ is
continuous with respect to local convergence on $\Pfp$, and similarly for
$\nucirc^\dgg_\de$.  Moreover, the function $\de\mapsto r(\La,a,\de)$ is
continuously differentiable on $(\de_{\rm c},\infty)$ and satisfies
\be\label{difr}
-\dif{\de}r(\La,a,\de)=
\frac{\nucirc_\de\captimes\nucirc_\de^\dgg\,(\{0\})}
{\lla\nucirc_\de\captimes\nucirc_\de^\dgg\rra}>0
\qquad\big(\de\in(\de_{\rm c},\infty)\big).
\ee
\et

\noi
\tbf{Remark} The continuity of $\nucirc_\de$ and $\nucirc^\dgg_\de$ as a
function of $\de$ in the sense of local convergence on $\Pfp$ is easily seen
to imply the continuity of the right-hand side of (\ref{difr}) in $\de$. On the
other hand, no such conclusion could be drawn from continuity in the sense of
vague convergence, since the functions $A\mapsto 1_{\{A=\{0\}\}}$ and
$A\mapsto |A|^{-1} 1_{\{0\in A\}}$ (which occur in the definition of
$\lla\,\cdot\,\rra$) are not continuous with respect to the
topology on $\Pc_+$.


\smallskip
The differentiability of the exponential growth rate in the subcritical regime
is expected. Indeed, for normal (unoriented) percolation in the subcritical
regime, it is even known that the number of open clusters per vertex and the
mean size of the cluster at the origin depend analytically on the percolation
parameter. This result is due to Kesten \cite{Kes81}; see also
\cite[Section~6.4]{Gri99}. For oriented percolation in one plus one dimension
in the \emph{supercritical} regime, Durrett \cite[Section~14]{Dur84} has shown
that the percolation probability is infinitely differentiable as a function of
the percolation parameter. It is not so clear, however, if the
methods in these papers can be adapted to cover the exponential growth
rate. At any rate, they would not give very explicit information about the
derivative such as positivity.

In principle, if for a given lattice one can show that the right-hand side of
(\ref{difr}) stays positive uniformly as $\de\down\de_{\rm c}$, then this
would imply that 
$r(\de)\sim(\de-\de_{\rm c})^1$ as $\de\down\de_{\rm c}$,
i.e., the critical exponent associated with the function $r$ is one. But this
is probably difficult in the most interesting cases, such as $\Z^d$ above the
critical dimension.

\subsection{Discussion and outlook}\label{S:discus}

This paper is part of a larger program, initiated in \cite{Swa09}, which aims
to describe all homogeneous eigenmeasures of $(\La,a,\de)$-contact processes
and to prove convergence for suitable starting measures. There are several
regimes of interest: the subcritical regime $\de>\de_{\rm c}$, the critical
regime $\de=\de_{\rm c}$, and the supercritical regime $\de<\de_{\rm c}.$ In
the supercritical regime one needs to distinguish further the case  $r=0$
(as for processes on $\Z^d$) and the case
$r>0$ (as for processes on trees).

In \cite{Swa09} some first, relatively weak results have been derived for
processes with $r=0$ in the supercritical regime. In particular, it was shown
that for such processes, there exists a unique homogeneous eigenmeasure with
eigenvalue zero \cite[Thm.~1.5]{Swa09}, but it has not been proved whether
there are homogeneous eigenmeasures with other eigenvalues, while convergence
has only been shown for one special initial measure and Laplace-transformed
times \cite[Corollary~3.4]{Swa09}.

Our present paper treats the subcritical case fairly conclusively. Arguably,
this should be the easiest regime. Indeed, our analysis is made easier by the
fact that the homogeneous eigenmeasures are concentrated on finite sets. As we
will see, such eigenmeasures are in one-to-one correspondence to
quasi-invariant laws for the contact process `modulo shifts'. More precisely,
call two sets $A,B\in\Pf$ equivalent if one is a translation of the other (see
(\ref{sim}) below), let $\ti A$ denote the corresponding equivalence class
containing $A$, and set $\tiPf:=\{\ti A:A\in\Pfp\}$. Then, for any
$(\La,a,\de)$-contact process $\eta$ started in a finite initial state,
$(\ti\eta_t)_{t\geq 0}$ is a Markov process with countable state space
$\tiPf$. In the subcritical regime, this process a.s.\ ends up in the trap
$\ti\emptyset$. We will show that there is a one-to-one correspondence between
eigenmeasures that are concentrated on $\Pfp$ and quasi-invariant laws for
$\ti\eta$. In particular, our results imply that the law of a subcritical
contact process modulo shifts, started in any finite initial state and
conditioned to be alive at time $t$, converges as $t\to\infty$ to a
quasi-invariant law (Theorem~\ref{T:quasi} below).

For certain discrete-time versions of the contact process on $\Z^d$, as well
as for some other, similar Markov chains, an analogous result has been proved
in \cite{FKM96}. Our methods differ significantly from the methods used there,
since we use eigenmeasures of the forward and dual process to construct
positive left and right eigenfunctions of the forward process. This simplifies
our proofs, but, since this approach makes essential use of contact process
duality, it is less generally applicable. The correspondence between
homogeneous eigenmeasures and quasi-invariant laws of the process modulo
shifts is only available in the subcritical regime. In contrast, in the
critical and supercritical regimes, we expect homogeneous eigenmeasures to be
concentrated on infinite sets, hence the techniques of the present paper are
not applicable.

Nevertheless, our methods give some hints on what to do in some of the other
regimes as well. Formula (\ref{difr}), which we expect to hold more generally,
says, roughly speaking, that $-\dif{\de}r(\La,a,\de)$ is the probability that
two independent sets, which are distributed according to the eigenmeasures
$\nucirc$ and $\nucirc^\dgg$ 
of the forward and dual (backward) process, and which are conditioned on
having nonempty intersection, intersect in a single point. In view of this, it
is tempting to try to replace the fact that $\nucirc$ and $\nucirc^\dgg$ are
each concentrated on finite sets, which
holds only in the subcritical regime, by the weaker assumption that the
intersection measure $\nucirc\captimes\nucirc^\dgg$ is concentrated on finite
sets. In particular, one wonders if this always holds in the regime $r>0$.

A simpler problem, which we have not pursued in the present paper, is to
investigate higher-order derivatives of $r(\La,a,\de)$ with respect to $\de$
or derivatives with respect to the infection rates $a(i,j)$. It seems likely
that the latter are strictly positive in the subcritical regime and given by a
formula similar to (\ref{difr}). Controlling higher-order derivatives of
$r(\La,a,\de)$ might be more difficult; in particular, we do not know if the
function $\de\mapsto r(\La,a,\de)$ is concave, or 
(which in view of (\ref{difr}) is a similar question), if the conditional laws
$\nucirc_\de(\,\cdot\,|\{A:0\in A\})$ are decreasing in the stochastic order,
as a function of~$\de$. 

\section{Main line of the proofs}

In this section we give an overview of the main line of our arguments.
In particular, we give the proofs of Theorems~\ref{T:eigcon}
and \ref{T:difr} in Sections~\ref{S:exunique} and
\ref{S:deriv} respectively. These proofs are based on a collection
of lemmas and propositions which are stated here but whose proofs are in
most cases postponed until later.

In short, the line of the arguments is as follows. We start in
Section~\ref{S:locfin2} by collecting some general facts about locally finite
measures on $\Pc_+$. In particular, we discuss the relation between vague and
local convergence, and we show that a homogeneous, locally finite measure on
$\Pfp$ can be seen as the `law' of a random finite set, shifted to a uniformly
chosen position in the lattice.

In Section~\ref{S:Existence}, we then prove the existence part of
Theorem~\ref{T:eigcon}. Since existence of an eigenmeasure with eigenvalue $r$
has already been proved in \cite{Swa09}, the main task is proving that 
there exists such an eigenmeasure that is moreover concentrated on
$\Pfp$. This is achieved by a covariance calculation.

Once existence is proved, we fix an eigenmeasure $\nucirc$ that is
concentrated on $\Pfp$, and likewise $\nucirc^{\dgg}$ for the dual process,
and set out to prove the convergence in (\ref{eeigcon}), which will then also
settle uniqueness. Our main strategy will be to divide out
translations and show that the resulting process modulo shifts is
$\la$-positive, which means that the subprobability kernels $P_t$ in
(\ref{Ptdef}) can be transformed, by a variation of Doob's $h$-transform, into
probability kernels belonging to a positively recurrent Markov process.
This part of the argument is carried out in Section~\ref{S:toquasi}.

To prepare for this, in Sections~\ref{S:tieta} and \ref{S:hmu}, we study left
and right eigenvectors of the semigroup of the $(\La,a,\de)$-contact process
modulo shifts. In particular, in Section~\ref{S:tieta}, we show that there is
a one-to-one correspondence between eigenmeasures that are concentrated on
finite sets and quasi-invariant laws (i.e., normalized positive left
eigenvectors) of the process modulo shifts. In Section~\ref{S:hmu}, we show
that moreover, each eigenmeasure of the dual $(\La,a^\dgg,\de)$-contact
process gives rise to a \emph{right} eigenvector for the $(\La,a,\de)$-contact
process modulo shifts.

In Section~\ref{S:exunique}, we use this to prove the convergence in
(\ref{eeigcon}), completing the proof of Theorem~\ref{T:eigcon}. We obtain
vague convergence for general starting measures by duality, using the
ergodicity of the Doob transform of the dual $(\La,a^\dgg,\de)$-contact
process modulo shifts. For starting measures that are concentrated on $\Pfp$,
we moreover obtain pointwise convergence by using the ergodicity of the Doob
transformed (forward) $(\La,a,\de)$-contact process modulo shifts, which
together with vague convergence, by a general lemma from
Section~\ref{S:locfin2}, implies local convergence on $\Pfp$.

In order to prove Theorem~\ref{T:difr}, in Section~\ref{S:cont} we show
continuity of the eigenmeasures $\nucirc$ in the recovery rate $\de$.
Continuity in the sense of vague convergence follows easily from a compactness
argument and uniqueness, but continuity in the sense of local convergence on
$\Pfp$ requires more work. We use a generalization of the covariance
calculation from Section~\ref{S:Existence} to obtain `local tightness', which
together with vague convergence, by a general lemma from
Section~\ref{S:locfin2}, implies local convergence on $\Pfp$.

In Section~\ref{S:deriv}, finally, we use the results proved so far to take the
limit $t\to\infty$ in Russo's formula (\ref{Russo}) and prove formula
(\ref{difr}), thereby completing the proof of Theorem~\ref{T:difr}.

At this point, the proofs of our main results are complete, but they depend on
a number of lemmas and propositions whose proofs have for readability been
postponed until later. We supply these in Section~\ref{S:detail}. The paper
concludes with two appendices. In Appendix~\ref{A:decay}, we point out how the
arguments in \cite{AJ07} generalize to the class of contact processes
considered in the present article.  Appendix~\ref{A:quasi} contains background
material on $\lambda$-positivity and quasi-invariant laws.

\subsection{More on locally finite measures}\label{S:locfin2}

In this section, we elaborate on the discussion in Section~\ref{S:locfin} of (contact processes started
in) locally finite measures on $\Pc_+$ by formulating
some lemmas that will be useful in what follows.

Recall from Section~\ref{S:locfin} the definition of vague convergence and of
local convergence on $\Pfp$, and recall that $\Pfi{i}:=\{A\in\Pf:i\in A\}$.
If $\mu_n,\mu$ are measures on $\Pfp$, then we say the $\mu_n$ converge to
$\mu$ \emph{pointwise on $\Pfp$} if $\mu_n(\{A\})\to\mu(\{A\})$ for all
$A\in\Pfp$. We say that the $(\mu_n)_{n\geq 1}$ are \emph{locally tight} if for
each $i\in\La$ and $\eps>0$ there exists a finite $\Di\sub\Pfi{i}$ such that
$\sup_n\mu_n(\Pfi{i}\beh\Di)\leq\eps$.  The next proposition, the proof of
which can be found in Section~\ref{S:lf}, connects all these definitions.

\bp\pbf{(Local convergence)}\label{P:locon}
Let $\mu_n,\mu$ be locally finite measures on $\Pc_+$ that are concentrated on
$\Pfp$. Then the following statements are equivalent.
\begin{enumerate}
\item $\mu_n\Rightarrow\mu$ locally on $\Pfp$.
\item $\mu_n\to\mu$ pointwise on $\Pfp$ and the $(\mu_n)_{n\geq 1}$ are
locally tight.
\item $\mu_n\Rightarrow\mu$ vaguely on $\Pc_+$ and the $(\mu_n)_{n\geq 1}$ are
locally tight.
\item $\mu_n\Rightarrow\mu$ vaguely on $\Pc_+$ and $\mu_n\to\mu$ pointwise on
$\Pfp$.
\end{enumerate}
\ep

Recall the definition of the intersection measure $\mu\captimes\nu$ in
(\ref{psidef}). The next lemma, whose proof can be found in
Section~\ref{S:lf}, says that the operation $\captimes$ is
continuous with respect to vague and local convergence.

\bl\pbf{(Intersection measure)}\label{L:intsect}
If $\mu,\nu$ are locally finite measures on $\Pc_+$, then $\mu\captimes\nu$ is
a locally finite measure on $\Pc_+$. If $\mu_n,\nu_n$ are locally finite
measures on $\Pc_+$ that converge vaguely to $\mu,\nu$, respectively, then
$\mu_n\captimes\nu_n$ converges vaguely to $\mu\captimes\nu$. If moreover
either the $\mu_n$ or the $\nu_n$ are concentrated on $\Pfp$ and converge
locally on $\Pfp$, then the $\mu_n\captimes\nu_n$ are concentrated on $\Pfp$
and converge locally on $\Pfp$.
\el

It is often useful to view a homogeneous, locally finite measure on $\Pfp$ as
the `law' of a random finite subset of $\La$, shifted to a uniformly chosen
position in $\La$. To formulate this precisely, we define an equivalence
relation on $\Pf$ by
\be\label{sim}
A\sim B\quad\mbox{iff}\quad A=iB\quad\mbox{for some }i\in\La,
\ee
and we let $\tiPf:=\{\ti A:A\in\Pf\}$ with $\ti A:=\{iA:i\in\La\}$ 
denote the set of equivalence classes. We can think of $\tiPf$ as the space of
finite subsets of the lattice `modulo shifts'. Recall the definition of
$\lla\mu\rra$ from (\ref{cmudef}). We have the following simple lemma, which
will be proved in Section~\ref{S:lf}.

\bl\pbf{(Homogeneous measures on the finite sets)}\label{L:homfin}
Let $\De$ be a $\Pfp$-valued random variable and let $c>0$. Then
\be\label{muDe}
\mu:=c\sum_{i\in\La}\P\big[i\De\in\cdot\,\big]
\ee
defines a nonzero, homogeneous measure on $\Pfp$ such that $\lla\mu\rra=c$.
The measure $\mu$ is locally finite if and only if $\E\big[|\De|\big]<\infty$.
Conversely, any nonzero, homogeneous measure on $\Pfp$ such that
$\lla\mu\rra<\infty$ can be written in the form (\ref{muDe}) with $c=\lla\mu\rra$ for some
$\Pfp$-valued random variable $\De$, and the law of $\ti\De$ is uniquely
determined by $\mu$.
\el

We finally turn our attention to contact processes started in infinite initial
`laws'. Recall the definition of the subprobability kernels $P_t$ in
(\ref{Ptdef}) and of the meaures $\mu P_t$ in (\ref{muPtdef}). We cite the
following simple fact from \cite[Lemma~3.3]{Swa09}.

\bl\pbf{(Process started in infinite law)}\label{L:infstart}
If $\mu$ is a homogeneous, locally finite measure on $\Pc_+$, then $\mu P_t$
is a homogeneous, locally finite measure on $\Pc_+$ for each $t\geq 0$.
If $\mu_n,\mu$ are homogeneous, locally finite measures on $\Pc_+$ such that
$\mu_n\Rightarrow\mu$, then $\mu_nP_t\Rightarrow\mu P_t$ for all $t\geq 0$,
where $\Rightarrow$ denotes vague convergence.
\el

\subsection{Existence of eigenmeasures concentrated on finite sets}
\label{S:Existence}

The first step in the proof of Theorem~\ref{T:eigcon} is to show that the
condition $r<0$ implies existence of a homogeneous eigenmeasure that is
concentrated on $\Pf$. 


We start by recalling how homogeneous eigenmeasures
with eigenvalue $r$ are constructed in \cite{Swa09}. For any
$(\La,a,\de)$-contact process, we can define homogeneous, locally finite
measures $\mu_t$ on $\Pc_+$ by
\be\label{mut}
\mu_t:=\sum_{i\in\La}\P[\eta^{\{i\}}_t\in\cdot\,]\big|_{\Pc_+}\qquad(t\geq 0).
\ee
We can think of $\mu_t$ as the law of a contact process started with one
infected site, distributed according to the counting measure on $\La$. It is
not hard to show (see \cite[formulas (3.8) and (3.20)]{Swa09}) that
\be\label{mupi}
\mu_t(\{A:0\in A\})=\E\big[|\eta^{\{0\}}_t|\big]=:\pi_t.
\ee
Let $\hat\mu_\la$ be the Laplace transform of $(\mu_t)_{t\geq 0}$, i.e.,
\be\label{Lap}
\hat\mu_\la:=\int_0^\infty\mu_t\:e^{-\la t}\di t\qquad(\la>r).
\ee
Then
\be\label{hatpi}
\hat\mu_\la(\{A:0\in A\})=\int_0^\infty\pi_t\:e^{-\la t}\di t=:\hat\pi_\la
\qquad(\la>r),
\ee
which is finite for $\la>r$ by the definition of the exponential growth rate
(see (\ref{rdef})). We cite the following result from
\cite[Corollary~3.4]{Swa09}, which yields the existence of homogeneous eigenmeasures.

\bp\pbf{(Convergence to eigenmeasure)}\label{P:conv}
The measures $\frac{1}{\hat\pi_\la}\hat\mu_\la$ $(\la>r)$ are relatively
compact in the topology of vague convergence of locally finite measures on
$\Pc_+$, and each subsequential limit as $\la\down r$ is a homogeneous
eigenmeasure of the $(\La,a,\de)$-contact process, with eigenvalue
$r(\La,a,\de)$.
\ep

We wish to show that for $r<0$, the approximation procedure in
Proposition~\ref{P:conv} yields an eigenmeasure that is concentrated on $\Pf$.
The key to this is the following lemma, which will be proved in
Section~\ref{S:covest} using a covariance calculation. Note that this lemma
still holds for general $r\in\R$.

\bl\pbf{(Uniform moment bound)}\label{L:unmom}
Let $\hat\mu_\la$ and $\hat\pi_\la$ be defined as in
(\ref{Lap})--(\ref{hatpi}).  Then, for any $(\La,a,\de)$-contact process with
exponential growth rate $r=r(\La,a,\de)$,
\be\label{umo}
\limsup_{\la\down r}\frac{1}{\hat\pi_\la}
\int\hat\mu_\la(\di A)1_{\{0\in A\}}|A|
\leq(|a|+\de)\int_0^\infty e^{-rt}\di t\,
\E\big[|\eta^{\{0\}}_t|\big]^2.
\ee
\el

As a consequence, we obtain the following result that completes the existence part of Theorem~\ref{T:eigcon}.

\bl\pbf{(Existence of an eigenmeasure on finite configurations)}\label{L:finex}
Assume that the exponential growth rate $r=r(\La,a,\de)$ of the
$(\La,a,\de)$-contact process satisfies $r<0$. Then there exists a homogeneous
eigenmeasure $\nucirc$ with eigenvalue $r$ of the $(\La,a,\de)$-contact process
such that
\be\label{nufin}
\int\nucirc(\di A)|A|1_{\{0\in A\}}<\infty.
\ee
\el
\tbf{Proof} By Proposition~\ref{P:conv}, we can choose $\la_n\down r$ such
that the measures $\frac{1}{\hat\pi_{\la_n}}\hat\mu_{\la_n}$ converge vaguely
to a homogeneous eigenmeasure $\nucirc$ with eigenvalue $r$. It follows from
(\ref{rdef}) that $\E[|\eta^{\{0\}}_t|]=e^{rt+o(t)}$ where $t\mapsto o(t)$ is
a continuous function such that $o(t)/t\to 0$ as $t\to\infty$, hence, by
(\ref{umo}), provided $r<0$,
\be\label{rneg}
\int_0^\infty e^{-rt}\di t\,\E\big[|\eta^{\{0\}}_t|\big]^2
=\int_0^\infty e^{2rt-rt+o(t)}\di t<\infty\qquad(r<0).
\ee
Let $\La_k$ be finite sets such that $0\in\La_k\sub\La$ and
$\La_k\up\La$. It is easy to check that $A\mapsto
f_k(A):=|A\cap\La_k|1_{\{0\in A\}}$ is a continuous, compactly supported real
function on $\Pc_+$. Therefore, by the vague convergence of
$\frac{1}{\hat\pi_{\la_n}}\hat\mu_{\la_n}$ to $\nucirc$, and by (\ref{umo}),
\be\ba{l}\label{numom}
\dis\int\nucirc(\di A)f_k(A)=\lim_{n\to\infty}
\frac{1}{\hat\pi_{\la_n}}\int\hat\mu_{\la_n}(\di A)f_k(A)\\[5pt]
\dis\quad\leq\liminf_{n\to\infty}
\frac{1}{\hat\pi_{\la_n}}\int\hat\mu_{\la_n}(\di A)|A|1_{\{0\in A\}}
\leq(|a|+\de)\int_0^\infty e^{-rt}\di t\,
\E\big[|\eta^{\{0\}}_t|\big]^2.
\ec
Letting $k\up\infty$, using the fact that the right-hand side
is finite by (\ref{rneg}), we arrive at (\ref{nufin}).\qed

\subsection{The process modulo shifts}\label{S:tieta}

Recall the definition of the equivalence relation in (\ref{sim}) and the
associated equivalence classes $\ti A$. Let $(\eta_t)_{t\geq 0}$ be a
$(\La,a,\de)$-contact process with $|\eta_0|<\infty$ a.s. Then
$\ti\eta=(\ti\eta_t)_{t\geq 0}$ is a Markov process with state space
$\tiPf$. We call $\ti\eta$ the \emph{$(\La,a,\de)$-contact process modulo
 shifts}. Clearly, the point $\ti\emptyset$ is a trap for this process. We
will prove that there is a one-to-one correspondence between eigenmeasures of
the $(\La,a,\de)$-contact process that are concentrated on $\Pfp$ and
quasi-invariant laws for the $(\La,a,\de)$-contact process modulo shifts.

To prepare for this, let $\ti P_t$ denote the transition probabilities of
$\ti\eta$, restricted to $\tiPfp$. Then $\ti P_t$ is a subprobability kernel
on $\tiPfp$ which is related to the subprobability kernel $P_t$ from
(\ref{Ptdef}) by
\be\label{PtiP}
\ti P_t(\ti A,\ti B)
=\!\!\sum\subb{C\in\Pfp}{\ti C=\ti B}\!\!P_t(A,C)
=m(B)^{-1}\sum_{i\in\La}P_t(A,iB)\qquad(t\geq 0,\ A,B\in\Pfp),
\ee
where we have to divide by the quantity
\be\label{mA}
m(A):=|\{i\in\La:iA=A\}|\qquad\big(A\in\Pfp\big)
\ee
to avoid double counting.\footnote{It is easy to see that the constant $m(A)$
 defined in (\ref{mA}) satisfies $m(A)\leq|A|$ and that $\{i\in\La:iA=A\}$
 is a finite subgroup of $\La$. If every element of $\La$ is of infinite
 order (as is the case, for example, for $\La=\Z^d$), then $m(A)=1$ for all
 finite $A\sub\La$.}

By Lemma~\ref{L:homfin}, each nonzero, homogeneous measure on $\Pfp$ such that
$\lla\mu\rra<\infty$ can be written as
\be
\mu=\lla\mu\rra\sum_{i\in\La}\P[i\De\in\cdot\,]
\ee
for some $\Pfp$-valued random variable $\De$. We write
\be\label{timu}
\ti\mu:=\P[\ti\De\in\cdot\,]
\ee
for the law of $\ti\De$. By Lemma~\ref{L:homfin}, $\ti\mu$ is uniquely
determined by $\mu$, and conversely, by (\ref{muDe}), $\ti\mu$ determines
$\mu$ up to a multiplicative constant.

We say that a function $f:\Pf\to\R$ is \emph{shift-invariant} if $f(iA)=f(A)$
for all $i\in\La$. For any shift-invariant function $f:\Pfp\to\R$, we let $\ti
f:\tiPfp\to\R$ denote the function defined by
\be\label{tif}
\ti f(\ti A):=f(A)\qquad(A\in\Pfp).
\ee
Then, clearly,
\be\label{tiPtif}
\ti P_t\ti f(\ti A)=\E[\ti f(\ti\eta^A_t)]=\E[f(\eta^A_t)]
=P_tf(A)\qquad(t\geq 0,\ A\in\Pfp).
\ee
The following simple lemma will be proved in Section~\ref{S:lf}.

\bl\pbf{(Laws on equivalence classes)}\label{L:timu}
Let $\mu$ be a nonzero, homogeneous measure on $\Pfp$ such that
$\lla\mu\rra<\infty$, let $\ti\mu$ be as in (\ref{timu}),
and let $P_t$ and $\ti P_t$ denote the transition probabilities of a
$(\La,a,\de)$-contact process and the latter modulo shifts, respectively.
Then
\be\label{mutimu}
\mu P_t(\{A\})=m(A)\lla\mu\rra\ti\mu\ti P_t(\ti A)\qquad(t\geq 0,\ A\in\Pfp).
\ee
Moreover, for any shift-invariant function $f:\Pfp\to\half$,
\be\label{timuf}
\sum_{\ti A\in\tiPfp}\ti\mu(\ti A)\ti f(\ti A)=\lla\mu\rra^{-1}\lla f\mu\rra.
\ee
\el

Formula (\ref{mutimu}) shows in particular that if $\mu$ is a homogeneous
measure on $\Pfp$ such that $\lla\mu\rra<\infty$ and $\la\in\R$, then
\be
\mu P_t=e^{\la t}\mu\quad(t\geq 0)\quad\mbox{if and only if}\quad
\ti\mu\ti P_t=e^{\la t}\ti\mu\quad(t\geq 0).
\ee
Here, the relation $\ti\mu\ti P_t=e^{\la t} \ti\mu$ says that the probability law
$\ti\mu$ is a \emph{quasi-invariant law} for the $(\La,a,\de)$-contact process
modulo shifts. Since $\mu$ determines $\ti\mu$ uniquely, and $\ti\mu$
determines $\mu$ uniquely up to a multiplcative constant, this shows that
there is a one-to-one correspondence between homogeneous eigenmeasures of the
$(\La,a,\de)$-contact process and quasi-invariant laws of the
$(\La,a,\de)$-contact process modulo shifts.

\subsection{Dual functions}\label{S:hmu}

We have just seen that eigenmeasures of the $(\La,a,\de)$-contact process give
rise to quasi-invariant laws for the $(\La,a,\de)$-contact process modulo
shifts, i.e., normalized, positive left eigenvectors of the semigroup $(\ti
P_t)_{t\geq 0}$. In the present section, we will see that moreover,
eigenmeasures of the \emph{dual} $(\La,a^\dgg,\de)$-contact process give rise
to positive \emph{right} eigenvectors of the semigroup $(\ti P_t)_{t\geq
 0}$. This will allow us to Doob-transform the $(\La,a,\de)$-contact process
modulo shifts into a positively recurrent Markov chain and use ergodicity of
the latter to prove, among other things, convergence to the quasi-invariant
law of the process conditioned not to have died out.

Recall the definition of the subprobability kernels $P_t$ and $P^\dgg_t$ in
(\ref{Ptdef}), which are the transition probabilities of the
$(\La,a,\de)$-contact process and its dual $(\La,a^\dgg,\de)$-contact process,
respectively, restricted to the set $\Pc_+$ of nonempty subsets of $\La$.
If $A\in\Pfp$, then $P_t(A,\,\cdot\,)$ is concentrated on the countable set
$\Pfp$, and we simply write $P_t(A,B):=P_t(A,\{B\})$. Let
\be\label{Gdom}
\Si(\Pfp):=\{f:\Pfp\to\R:
|f(A)|\leq K|A|^k+M\mbox{ for some }K,M,k\geq 0\}.
\ee
denote the class of real functions on $\Pfp$ of polynomial growth.
For any $f\in\Si(\Pfp)$, we define
\be
P_tf(A):=\sum_{B\in\Pfp}P_t(A,B)f(B)\qquad(t\geq 0,\ A\in\Pfp).
\ee
Then \cite[Prop.~2.1]{Swa09} implies that $P_t$ maps the space $\Si(\Pfp)$
into itself.

For any locally finite measure $\mu$ on $\Pc_+$, we define a function
$h_\mu:\Pfp\to\half$ by
\be\label{hdef}
h_\mu(A):=\int_{\Pc_+}\!\mu(\di B)\,1_{\{A\cap B\neq\emptyset\}}
\qquad(A\in\Pfp),
\ee
which is finite by Lemma~\ref{L:locfin}. We say that a function $f:\Pf\to\R$
is \emph{monotone} if $A\sub B$ implies $f(A)\leq f(B)$, and
\emph{subadditive} if $f(A\cup B)\leq f(A)+f(B)$, for all $A,B\in\Pf$. Below,
$\mu P_t$ is defined as in (\ref{muPtdef}).

\bl\pbf{(Linear bounds)}\label{L:hlin}
If $\mu$ be a nonzero, homogeneous, locally finite measure $\mu$ on $\Pc_+$.
For each nonzero, homogeneous, locally finite measure $\mu$ on $\Pc_+$, the
function $h_\mu$ in (\ref{hdef}) is shift-invariant, monotone, subadditive,
strictly positive on $\Pfp$, and satisfies $h_\mu(A)\leq h(\{0\})|A|$
$(A\in\Pfp)$. If $\mu$ is moreover concentrated on $\Pfp$, then
\be\label{hlin}
\lla\mu\rra|A|\leq h_\mu(A)\leq h(\{0\})|A|
\qquad(A\in\Pfp).
\ee
\el
\tbf{Proof} Shift-invariance, monotonicity, subadditivity and positivity are
easy to check; see \cite[Lemma~3.5]{Swa09}. Subadditivity and shift-invariance
now imply the upper bound in (\ref{hlin}). If $\mu$ is concentrated on $\Pfp$,
then by Lemma~\ref{L:homfin}, there exists a $\Pfp$-valued random variable
$\De$ such that $\mu$ can be written as in (\ref{muDe}). Letting $\kappa$ be a
$\La$-valued random variable such that $\kappa\in\De$ a.s., we observe that
\be
\dis h_\mu(A)
=\lla\mu\rra\sum_{i\in\La}\P[A\cap i\De\neq\emptyset]
\geq\lla\mu\rra\sum_{i\in\La}\P[A\cap\{i\kappa\}\neq\emptyset]
=\lla\mu\rra|A|.
\ee
\qed

The next lemma is a simple consequence of duality.

\bl\pbf{(Dual function)}\label{L:dufun}
For each nonzero, homogeneous, locally finite measure $\mu$ on $\Pc_+$,
one has
\be\label{hmuP}
h_{\mu P_t}=P^\dgg_th_\mu\qquad(t\geq 0).
\ee
\el
\tbf{Proof} The upper bound in (\ref{hlin}) shows that $h_\mu\in\Si(\Pfp)$, so
$P^\dgg_th_\mu$ is well-defined. Also, $\mu P_t$ is homogeneous and locally
finite by Lemma~\ref{L:infstart}, so that $h_{\mu P_t}$ is well-defined. Now
\be\ba{l}
\dis h_{\mu P_t}(A)=\int_{\Pc_+}\!\mu P_t(\di B)\,1_{\{A\cap B\neq\emptyset\}}
=\int_{\Pc_+}\!\mu(\di C)\P[A\cap\eta^C\neq\emptyset]\\[5pt]
\dis\quad=\int_{\Pc_+}\!\mu(\di C)\P[\eta^{\dgg\,A}_t\cap C\neq\emptyset]
=\!\sum_{B\in\Pfp}\!P^\dgg_t(A,B)\int_{\Pc_+}\!\mu(\di C)\,1_{\{A\cap C\neq\emptyset\}}
=P^\dgg_th_\mu(A).
\ec
\qed

In particular, if $\mu$ is a homogeneous eigenmeasure of the
$(\La,a^\dgg,\de)$-contact process with eigenvalue $\la$, then
Lemma~\ref{L:dufun} implies that $P_th_\mu=e^{\la t}h_\mu$, i.e., $h_\mu$ is a
right eigenfunction of $P_t$, for each $t\geq 0$. By formula (\ref{tiPtif}),
$h_\mu$ then also gives rise to a right eigenfunction $\ti h_\mu$ of the
semigroup of the process modulo shifts. The following lemma, whose proof 
can be found in Section~\ref{S:lf}, will be handy in what follows.

\bl\pbf{(Intersection and weighted measures)}\label{L:cpsi}
Let $\mu,\nu$ be homogeneous locally finite measures on $\Pc_+$, assume that
$\mu$ is concentrated on $\Pfp$, and let $h_\nu$ be defined as in
(\ref{hdef}). Then
\be\label{cpsi}
\lla\mu\captimes\nu\rra=\lla h_\nu\mu\rra.
\ee
If moreover $\int\mu(\di A)|A|1_{\{0\in A\}}<\infty$, then $h_\nu\mu$ is
locally finite.
\el

\subsection{Convergence to the quasi-invariant law}\label{S:toquasi}

By Theorem~\ref{T:rprop}~(a), the $(\La,a,\de)$-contact processes and its dual
$(\La,a^\dgg,\de)$-contact processes have the same exponential growth rate
$r=r(\La,a,\de)=r(\La,a^\dgg,\de)$. In particular, if $r<0$, then by
Lemma~\ref{L:finex}, there exist homogeneous eigenmeasures $\nucirc$ and
$\nucirc^\dgg$ of the $(\La,a,\de)$- and $(\La,a^\dgg,\de)$-contact process,
respectively, both with eigenvalue $r$, such that
\be\label{nunufin}
\int\nucirc(\di A)|A|1_{\{0\in A\}}<\infty
\quad\mbox{and}\quad 
\int\nucirc^\dgg(\di A)|A|1_{\{0\in A\}}<\infty.
\ee 
We normalize $\nucirc$ and $\nucirc^\dgg$ such that $\int\nucirc(\di
A)1_{\{0\in A\}}=1=\int\nucirc^\dgg(\di A)1_{\{0\in A\}}$. For the moment, we
do not know yet if $\nucirc$ and $\nucirc^\dgg$ are unique,
but we simply fix any two such measures and we let
\be\label{tinu}
\ti\nu:=\ti\nucirc\quand\ti\nu^\dgg:=\ti\nucirc^\dgg
\ee
denote the associated probability laws on $\tiPfp$ as in (\ref{timu}).  We let
$h_\nucirc$ and $h_{\nucirc^\dgg}$ denote the associated dual functions as in
(\ref{hdef}) and let $\ti h_\nucirc$ and $\ti h_{\nucirc^\dgg}$ denote the
associated functions on $\tiPfp$ as in (\ref{tif}).  Finally, we define
\emph{Doob-transformed}\footnote{Doob's classical $h$-transform is based on a
 positive harmonic function $h$. In (\ref{QQ}) we use a slight generalization
 of this where $h$ is a positive eigenfunction of the generator. This is a
 special case of what is called a `compensated $h$-transform' in
 \cite[Lemma~3]{FS02}.} probability kernels on $\tiPfp$ by
\be\label{QQ}
Q_t(\ti A,\ti B):=
e^{-rt}\ti h_{\nucirc^\dgg}(\ti A)^{-1}P_t(\ti A,\ti B)\ti h_{\nucirc^\dgg}(\ti B)
\qquad(t\geq 0,\ \ti A,\ti B\in\tiPfp).
\ee
The results from the last two sections, together with classical results about
quasi-invariant laws, then combine to give the following result.

\begin{atheorem}\pbf{(Convergence to the quasi-invariant law)}
\label{T:quasi}
Assume that the infection rates sa\-tisfy the irreducibility condition
(\ref{irr}) and that the exponential growth rate from (\ref{rdef}) satisfies
$r<0$. Then:
\begin{itemize}
\item[\rm(a)] The $(Q_t)_{t\geq 0}$ from (\ref{QQ}) are the transition
probabilities of a positively recurrent con\-ti\-nu\-ous-time Markov chain
with unique invariant law $\ti\pi$ given by
\be\label{pip}
\ti\pi(\ti A):=\frac{\lla\nucirc\rra}{\lla\nucirc\captimes\nucirc^\dgg\rra}
\ti\nu(\ti A)\ti h_{\nucirc^\dgg}(\ti A)\qquad(\ti A\in\tiPfp).
\ee
\item[\rm(b)] The law $\ti\nu$ is a quasi-invariant law for the
$(\La,a,\de)$-contact process modulo shifts. Moreover, for any $A\in\Pfp$, one
has
\be
\P[\ti\eta^A_t\in\cdot\,|\,\eta^A_t\neq\emptyset]\Asto{t}\ti\nu,
\ee
where $\Rightarrow$ denotes weak convergence of probability laws on $\tiPfp$.
\end{itemize}
\end{atheorem}
\tbf{Proof} By Lemmas~\ref{L:timu} and \ref{L:dufun}, the measure $\ti\nu$ and
function $\ti h_{\nucirc^\dgg}$ are left and right eigenvectors of the
operators $\ti P_t$, i.e.,
\be\label{leftright}
\ti\nu\ti P_t=e^{rt}\ti\nu\quand\ti P_t\ti h_{\nucirc^\dgg}=e^{rt}\ti h_{\nucirc^\dgg}
\qquad(t\geq 0).
\ee
Moreover, by formula (\ref{timuf}) and Lemma~\ref{L:cpsi},
\be\label{summable}
\sum_{\ti A\in\tiPfp}\ti\nu(\ti A)\ti h_{\nucirc^\dgg}(\ti A)
=\lla\nucirc\rra^{-1}\lla h_{\nucirc^\dgg}\nucirc\rra
=\lla\nucirc\rra^{-1}\lla\nucirc\captimes\nucirc^\dgg\rra<\infty.
\ee
Thus, we have found positive left and right eigenfunctions of $(\ti
P_t)_{t\geq 0}$ whose pointwise product is summable. The statements of the
theorem now follow readily by well-known methods. More precisely, parts~(a)
and (b) follow from Lemmas~\ref{L:Rpos} and \ref{L:quasi} in the appendix,
respectively. As explained in the proof of Lemma~\ref{L:Rpos} there, formulas
(\ref{leftright}) and (\ref{summable}) imply that the $(\La,a,\de)$-contact
process modulo shifts is $\la$-positive in the sense of Kingman
\cite{Kin63}. Note that the use of $\la$-positivity and its discrete time
analogue R-positivity in the study of quasi-invariant laws is well-known, see
e.g.\ \cite{FKM96}.\qed

\subsection{Convergence to the eigenmeasure}\label{S:exunique}

\tbf{Proof of Theorem~\ref{T:eigcon}} The existence of $\nucirc$ and
$\nucirc^\dgg$ has already been proved in Lemma~\ref{L:finex}, so uniqueness
will follow once we prove the convergence in (\ref{eeigcon}), with the
$\nucirc$ that we fixed earlier. 
We need to prove two statements: vague convergence for general (nonzero,
homogeneous, locally finite) initial measures $\mu$ and local convergence on
$\Pfp$ if $\mu$ is concentrated on $\Pfp$.

We start with vague convergence. By Lemma~\ref{L:locfin}, it suffices to show
that
\be\label{conv1}
e^{-rt}\int\mu P_t(\di A)1_{\{A\cap B\neq\emptyset\}}\asto{n}
c\int\nucirc(\di A)1_{\{A\cap B\neq\emptyset\}}
\qquad(B\in\Pfp),
\ee
where $c>0$ is given in (\ref{cform}). By Lemma~\ref{L:dufun}, we
observe that
\be\ba{l}\label{conv2}
\dis e^{-rt}\int\mu P_t(\di A)1_{\{A\cap B\neq\emptyset\}}
=e^{-rt}h_{\mu P_t}(B)=e^{-rt}P^\dgg_th_\mu(B)\\[5pt]
=\E[h_\mu(\eta^{\dgg\,B}_t)]
=\E[\ti h_\mu(\ti\eta^{\dgg\,B}_t)]=\ti P^\dgg_t\ti h_\mu(\ti B)
\qquad(B\in\Pfp),
\ec
where $(\ti P^\dgg_t)_{t\geq 0}$ are defined as in (\ref{PtiP}) but for the
$(\La,a^\dgg,\de)$-contact process modulo shifts. Applying
Theorem~\ref{T:quasi}~(a) to the dual process, we obtain that
\be
Q^\dgg_t(\ti A,\ti B):=
e^{-rt}\ti h_\nucirc(\ti A)^{-1}P^\dgg_t(\ti A,\ti B)\ti h_\nucirc(\ti B)
\qquad(t\geq 0,\ \ti A,\ti B\in\tiPfp)
\ee
are the transition probabilities of an irreducible, positively recurrent
Markov process with state space $\tiPfp$ and invariant law $\pi^\dgg$ given by
\be\label{pipi}
\ti\pi^\dgg(\ti A)
:=\frac{\lla\nucirc^\dgg\rra}{\lla\nucirc\captimes\nucirc^\dgg\rra}
\ti\nu^\dgg(\ti A)\ti h_\nucirc(\ti A)\qquad(\ti A\in\tiPfp).
\ee
The right-hand side of (\ref{conv2}) can now be rewritten as
\be
\ti P^\dgg_t\ti h_\mu(\ti B)
=\ti h_\nucirc(\ti B)\!\sum_{\ti A\in\tiPfp}\!Q^\dgg_t(\ti B,\ti A)
\frac{\ti h_\mu(\ti A)}{\ti h_\nucirc(\ti A)}
\qquad(\ti A\in\tiPfp).
\ee
By formula (\ref{hlin}) from Lemma~\ref{L:hlin}, $\ti h_\mu/\ti h_\nucirc$ is
a bounded function, so we may use the ergodicity of the irreducible,
positively recurrent Markov process with transition probabilities
$(Q^\dgg_t)_{t\geq 0}$ to conclude that
\be\ba{l}\label{ticon}
\dis\ti P^\dgg_t\ti h_\mu(\ti B)
\asto{t}\ti h_\nucirc(\ti B)\!\sum_{\ti A\in\tiPfp}\!
\ti\pi^\dgg(\ti A)\frac{\ti h_\mu(\ti A)}{\ti h_\nucirc(\ti A)}\\[5pt]
\dis\quad=\ti h_\nucirc(\ti B)
\frac{\lla\nucirc^\dgg\rra}{\lla\nucirc\captimes\nucirc^\dgg\rra}
\!\sum_{\ti A\in\tiPfp}\!\ti\nu^\dgg(\ti A)\ti h_\mu(\ti A)
\,=\,
\frac{\lla\mu\captimes\nucirc^\dgg\rra}{\lla\nucirc\captimes\nucirc^\dgg\rra}
\int\nucirc(\di A)1_{\{A\cap B\neq\emptyset\}},
\ec
where in the last step we have used formula (\ref{timuf}) from
Lemma~\ref{L:timu} as well as Lemma~\ref{L:cpsi}. Combining (\ref{ticon}) with
(\ref{conv1})--(\ref{conv2}), this proves the vague convergence in
(\ref{eeigcon}).

It remains to show that vague convergence can be
strengthened to local convergence on $\Pfp$ if $\mu$ is concentrated on
$\Pfp$. By Proposition~\ref{P:locon}~(iv), it suffices to prove pointwise
convergence. Let $Q_t$ and $\pi$ be given by
(\ref{QQ})--(\ref{pip}). Then formula (\ref{mutimu}) from
Lemma~\ref{L:timu} tells us that
\be
\mu P_t(\{B\})=m(B)\lla\mu\rra\ti\mu\ti P_t(\ti B)
=m(B)\lla\mu\rra\, e^{rt}\!\!\sum_{\ti A\in\tiPfp}\!\ti\mu(\ti A)
\ti h_{\nucirc^\dgg}(\ti A)Q_t(\ti A,\ti B)\ti h_{\nucirc^\dgg}(\ti B)^{-1}.
\ee
Here, on the right-hand side, we evolve the measure $\ti
h_{\nucirc^\dgg}\ti\mu$ under the semigroup $Q_t$. Using formula (\ref{timuf}) from
Lemma~\ref{L:timu} and Lemma~\ref{L:cpsi} we see that this measure is
finite with total mass given by
\be
\!\sum_{\ti A\in\tiPfp}\!\ti\mu(\ti A)\ti h_{\nucirc^\dgg}(\ti A)
=\lla\mu\rra^{-1}\lla h_{\nucirc^\dgg}\mu\rra
=\lla\mu\rra^{-1}\lla\mu\captimes\nucirc^\dgg\rra<\infty.
\ee
Using this and the ergodicity of the Markov process with transition
probabilities $(Q_t)_{t\geq 0}$, we find that
\be\ba{l}
\dis e^{-rt}\mu P_t(\{B\})\asto{t}m(B)\lla\mu\captimes\nucirc^\dgg\rra
\pi(\ti B)\ti h_{\nucirc^\dgg}(\ti B)^{-1}\\[5pt]
\dis\quad=
\frac{\lla\mu\captimes\nucirc^\dgg\rra}{\lla\nucirc\captimes\nucirc^\dgg\rra}
m(B)\lla\nucirc\rra\ti\nu(\ti B)
=\frac{\lla\mu\captimes\nucirc^\dgg\rra}{\lla\nucirc\captimes\nucirc^\dgg\rra}
\nucirc(\{B\}),
\ec
where in the last step we have used formula (\ref{mutimu}) from
Lemma~\ref{L:timu} (with $t=0$).\qed

\subsection{Continuity in the recovery rate}\label{S:cont}

The first step in proving Theorem~\ref{T:difr} will be to show continuity of
the map $(\de_{\rm c},\infty)\ni\de\mapsto\nucirc_\de$. We start by proving
continuity with respect to vague convergence, which is based on the following
abstract result, whose proof can be found in Section~\ref{S:infstart}.

\bl\pbf{(Limits of eigenmeasures)}\label{L:eiglim}
Let $\nu_n$ $(n\geq 0)$ be homogeneous eigenmeasures of
$(\La,a,\de_n)$-contact processes, with eigenvalues $\la_n$, normalized such
that $\int\nu_n(\di A)1_{\{0\in A\}}=1$. Assume that $\la_n\to\la$ and
$\de_n\to\de$. Then the $(\nu_n)_{n\geq 0}$ are relatively compact in the
topology of vague convergence, and each vague cluster point $\nu$ is a
homogeneous eigenmeasure of the $(\La,a,\de)$-contact processes, with
eigenvalue $\la$.
\el

Continuity of the map $(\de_{\rm c},\infty)\ni\de\mapsto\nucirc_\de$ is now a
simple consequence of Theorem~\ref{T:eigcon} and Lemma~\ref{L:eiglim}.

\bp\pbf{(Vague continuity of the eigenmeasure)}\label{P:eigcont}
Assume that the infection rates satisfy the irreducibility condition
(\ref{irr}). For $\de\in(\de_{\rm c},\infty)$, let $\nucirc_\de$ denote the
unique homogeneous eigenmeasure of the $(\La,a,\de)$-contact process
normalized such that $\int\nucirc_\de(\di A)1_{\{0\in A\}}=1$. Then the map
$\de\mapsto\nucirc_\de$ is continuous on $(\de_{\rm c},\infty)$ w.r.t.\ vague
convergence of locally finite measures on $\Pc_+$.
\ep
\tbf{Proof} Choose $\de_n,\de\in(\de_{\rm c},\infty)$ such that
$\de_n\to\de$. Since the eigenvalue $r(\La,a,\de)$ of the homogeneous
eigenmeasure $\nucirc_\de$ is continuous in $\de$ by
Theorem~\ref{T:rprop}~(b), Lemma~\ref{L:eiglim} implies that the measures
$(\nucirc_{\de_n})_{n\geq 0}$ are relatively compact in the topology of vague
convergence, and each vague cluster point is a homogeneous eigenmeasure of the
$(\La,a,\de)$-contact processes with eigenvalue $r(\La,a,\de)$. By
Theorem~\ref{T:eigcon}, this implies that $\nucirc_\de$ is the only vague
cluster point, hence the $\nucirc_{\de_n}$ converge vaguely to
$\nucirc_\de$.\qed

Unfortunately, continuity with respect to vague convergence is not enough to
prove continuity of the right-hand side of (\ref{difr}), and hence of the
derivative $\dif{\de}r(\La,a,\de)$. As mentioned earlier, we will remedy this
by proving continuity of the map $(\de_{\rm
c},\infty)\ni\de\mapsto\nucirc_\de$ with respect to local convergence on
$\Pfp$. Since vague convergence is already proved, by
Proposition~\ref{P:locon}~(iii), it suffices to prove local tightness.  This
is the most technical part of our proofs, since it involves estimating how
`large' the finite sets can be that $\nucirc_\de$ is concentrated on.  The
first step is to introduce a suitable concept of distance.  The next result
will be proved in Section~\ref{S:expmom}.

\bl\pbf{(Slowly growing metric)}\label{L:slowd}
Let $\La$ be a countable group and let $a:\La\times\La\to\half$ satisfy
(\ref{assum}). Then there exists a metric $d$ on $\La$ such that
\be\ba{rll}\label{dprop}
{\rm(i)}&\dis d(i,j)=d(ki,kj)&\dis\qquad(i,j,k\in\La),\\[5pt]
{\rm(ii)}&\dis\big|\{i\in\La:d(0,i)\leq M\}|<\infty
&\dis\qquad(0\leq M<\infty),\\[5pt]
{\rm(iii)}&\dis K_\ga(\La,a):=\sum_ia(0,i)e^{\ga d(0,i)}<\infty
&\dis\qquad(0\leq\ga<\infty).
\ec
\el

Next, we fix a metric $d$ as in (\ref{dprop}) and for each $0\leq\ga<\infty$,
we define a function $e_\ga:\Pf\to\half$ by
\be\label{ega}
e_\ga(A):=\sum_{i\in A}e^{\ga d(0,i)}\qquad(\ga\geq 0,\ A\in\Pf).
\ee
We note that a similar (but not entirely identical) function has proved useful
in the study of contact processes on trees, see \cite[formula~(I.4.3)]{Lig99}.
We have in particular $e_0(A)=|A|$. The next lemma says that there is a
well-defined exponential growth rate $r_\ga(\La,a,\de)$ associated with the
function $e_\ga$, which converges to our well-known exponential growth rate
$r(\La,a,\de)$ as $\ga\down 0$. The proof can be found in
Section~\ref{S:expmom}.

\bl\hspace*{-3pt}\pbf{(Exponential growth rates)}\hspace*{-2pt}\label{L:rga}
Let $(\eta^{\{0\}}_t)_{t\geq 0}$ be the $(\La,a,\de)$-contact process started
in $\eta^{\{0\}}_0=\{0\}.$ Let $d$ be a metric on $\La$ as in
Lemma~\ref{L:slowd}, and let $e_\ga$ be the function defined in (\ref{ega}).
Then, for each $0\leq\ga<\infty$, the limit
\be\label{rga}
r_\ga=r_\ga(\La,a,\de)
:=\lim_{t\to\infty}\ffrac{1}{t}\log\E\big[e_\ga(\eta^{\{0\}}_t)\big]
=\inf_{t>0}\ffrac{1}{t}\log\E\big[e_\ga(\eta^{\{0\}}_t)\big]
\ee
exists. The function $\ga\mapsto r_\ga$ is nondecreasing, right-continuous,
and satisfies
\be
-\de\leq r_\ga(\La,a,\de)\leq K_\ga(\La,a)\qquad(\ga\geq 0),
\ee 
where $K_\ga(\La,a)$ is defined in (\ref{dprop}).
\el

We can generalize the proof of Lemma~\ref{L:unmom} to yield a more general
version of that lemma (see Lemma~\ref{L:unexmom} below), which after taking
the limit (as in (\ref{numom})) yields the following bound on the
eigenmeasures $\nucirc_\de$. (We refer to Section~\ref{S:covest} for the
detailed proof.)

\bl\pbf{(Tightness estimate)}\label{L:nutight}
Let $(\eta^{\{0\}}_t)_{t\geq 0}$ be the $(\La,a,\de)$-contact process started
in $\eta^{\{0\}}_0=\{0\}$, let $r(\de)=r(\La,a,\de)$ be its exponential growth
rate, let $d$ be a metric on $\La$ as in Lemma~\ref{L:slowd}, and let $e_\ga$
be the function defined in (\ref{ega}). For $\de\in(\de_{\rm c},\infty)$, let
$\nucirc_\de$ denote the unique homogeneous eigenmeasure of the
$(\La,a,\de)$-contact process normalized such that $\int\nucirc_\de(\di
A)1_{\{0\in A\}}=1$. Then
\be\label{nutight}
\int\nucirc_\de(\di A)1_{\{0\in A\}}e_\ga(A)
\leq(|a|+\de)\int_0^\infty e^{-r(\de)t}\di t
\,\E\big[e_\ga(\eta^{\de,\,\{0\}}_t)\big]^2
\qquad\big(\ga\geq 0,\ \de\in(\de_{\rm c},\infty)\big).
\ee
\el
\noi
With this preparation we are now ready to prove the desired local continuity.
\bp\pbf{(Local continuity of the eigenmeasure)}\label{P:eigpoint}
Assume that the infection rates satisfy the irreducibility condition
(\ref{irr}). For $\de\in(\de_{\rm c},\infty)$, let $\nucirc_\de$ denote the
unique homogeneous eigenmeasure of the $(\La,a,\de)$-contact process
normalized such that $\int\nucirc_\de(\di A)1_{\{0\in A\}}=1$. Then the map
$\de\mapsto\nucirc_\de$ is continuous on $(\de_{\rm c},\infty)$ in the sense
of local convergence on $\Pfp$.
\ep
\tbf{Proof} Vague continuity of the map
$(\de_{\rm c},\infty)\ni\de\mapsto\nucirc_\de$ has been proved in
Proposition~\ref{P:eigcont}, so by Proposition~\ref{P:locon}~(iii), it
suffices to show that for any $\de_\ast\in(\de_{\rm c},\infty)$ there exists
an $\eps>0$ such that the measures
$(\nucirc_\de)_{\de\in(\de_\ast-\eps,\de_\ast+\eps)}$ are locally tight.

By property~(\ref{dprop})~(ii), for each $\ga>0$ and $K<\infty$, the set
$\{A\in\Pfi{0}:e_\ga(A)\leq K\}$ is finite. Thus, by Lemma~\ref{L:nutight}, to
prove the required local tightness, it suffices to show that for each
$\de_\ast\in(\de_{\rm c},\infty)$ there exist a $\ga>0$ and $\eps>0$ such that
\be\label{unint}
\sup_{\de\in(\de_\ast-\eps,\de_\ast+\eps)}
\int_0^\infty e^{-r(\de)t}\di t\,\E\big[e_\ga(\eta^{\de,\,\{0\}}_t)\big]^2
<\infty.
\ee
By the continuity of $\de\mapsto r(\de)$ (Theorem~\ref{T:rprop}~(b)), we can
choose $\eps>0$ such that $\de_{\rm c}<\de_\ast-\eps$ and
\be
r(\de_\ast-\eps)\leq\frac{4}{5}r(\de_\ast+\eps).
\ee
Let $r_\ga=r_\ga(\de)$ be the exponential growth rate associated with the
function $e_\ga$. By Lemma~\ref{L:rga}, the function $\ga\mapsto r_\ga$ is
right-continuous, so we can choose $\ga>0$ such that
\be
r_\ga(\de_\ast-\eps)\leq\frac{3}{4}r(\de_\ast-\eps).
\ee
By the fact that $r(\de)$ is nonincreasing in $\de$ and the law of
$\eta^{\de,\,\{0\}}_t$ is nonincreasing in $\de$ with respect to the
stochastic order, it follows that for all
$\de\in(\de_\ast-\eps,\de_\ast+\eps)$,
\be\ba{l}
\dis\int_0^\infty e^{-r(\de)t}\di t\,\E\big[e_\ga(\eta^{\de,\,\{0\}}_t)\big]^2
\leq\int_0^\infty e^{-r(\de_\ast+\eps)t}\di t\,
\E\big[e_\ga(\eta^{\de_\ast-\eps,\,\{0\}}_t)\big]^2\\[5pt]
\dis\quad=\int_0^\infty \di t\,
\ex{(2r_\ga(\de_\ast-\eps)-r(\de_\ast+\eps))t+o(t)}
\leq\int_0^\infty \di t\,
\ex{\frac{1}{5}r(\de_\ast+\eps)t+o(t)}<\infty,
\ec
where $t\mapsto o(t)$ is continuous, $o(t)/t\to 0$ for $t\to\infty$ by the
definition of $r_\ga$ in Lemma~\ref{L:rga}, and we have used that
$2r_\ga(\de_\ast-\eps)\leq 2\cdot\frac{3}{4}\cdot\frac{4}{5}r(\de_\ast+\eps)
=\frac{6}{5}r(\de_\ast+\eps)$. This proves (\ref{unint}) and hence the required
local tightness.\qed

\subsection{The derivative of the exponential growth rate}\label{S:deriv}

Let us define homogeneous, locally finite measures $\chi_A$ on
$\Pfp$ by
\be\label{chidef}
\chi_A:=\sum_{i\in\La}\de_{iA}\qquad(A\in\Pfp),
\ee
where $\de_{iA}$ denotes the delta measure on $\Pfp$ at the point $iA$.
Let $(P^\de_t)_{t\geq 0}$ and $(P^{\dgg\,\de}_t)_{t\geq 0}$ be the
subprobability kernels defined in (\ref{Ptdef}) for the $(\La,a,\de)$- and
$(\La,a^\dgg,\de)$-contact processes, respectively, in dependence on $\de$.
Note that $\chi_{\{0\}}P^\de_t$ denotes the `law' at time $t$ of the process
started with a single infected site distributed according to the counting
measure on $\La$. We start by rewriting Russo's formula (\ref{Russo}) in
terms of the objects we are working with.

\bl\pbf{(Differential formula)}
For each $t\geq 0$, the function $\half\ni\de\mapsto\E[|\eta^{\de,\{0\}}_t|]$
is continuously differentiable and satisfies
\be\label{Russo2}
-\dif{\de}\frac{1}{t}\log\E\big[|\eta^{\de,\,\{0\}}_t|\big]
=\frac{1}{t}\int_0^t\!\di s\,\frac{\dis
\chi_{\{0\}}P^\de_s\captimes\chi_{\{0\}} P^{\dgg\,\de}_{t-s}\,(\{0\})}
{\dis\lla\chi_{\{0\}}P^\de_s\captimes
\chi_{\{0\}} P^{\dgg\,\de}_{t-s}\rra}.
\ee
\el
\tbf{Proof} By (\ref{Russo}) and the definition of the Campbell law $\hat\P_t$
in (\ref{Campbell})
\be
\dis-\dif{\de}\frac{1}{t}\log\E\big[|\eta^{\de,\,\{0\}}_t|\big]
=\frac{1}{t}
\int_0^t\di s\;\frac{1}{\E\big[|\eta^{\de,\,\{0\}}_t|\big]}
\sum_{i,j}\P[(0,0)\leadsto_{(j,s)}(i,t)],
\ee
where
\be\ba{l}
\dis\sum_{i,j}\P[(0,0)\leadsto_{(j,s)}(i,t)]
=\sum_{i,j}\P[(j^{-1},-s)\leadsto_{(0,0)}(j^{-1}i,t-s)]\\[5pt]
\dis\quad=\sum_{i,j}
\P[\eta^{\de\,\{i\}}_s\cap\eta^{\dgg\,\de\,\{j\}}_{t-s}=\{0\}]
=\int\!\chi_{\{0\}}P^\de_s\,(\di A)\,
\int\!\chi_{\{0\}}P^{\dgg\,\de}_{t-s}\,(\di B)\,1_{\{A\cap B=\{0\}\}},
\ec
and for $0 \leq s \leq t,$
\be\ba{l}
\dis\E\big[|\eta^{\de,\,\{0\}}_t|\big]
=\sum_i\P\big[\eta^{\de\,\{0\}}_t\cap\{i\}\neq\emptyset\big]
=\sum_i\P\big[\eta^{\de\,\{i^{-1}\}}_t\cap\{0\}\neq\emptyset\big]\\[5pt]
\dis\quad=\sum_{i,j}
\E\big[|\eta^{\de\,\{i^{-1}\}}_t\cap\{j\}|^{-1}
1_{\{0\in\eta^{\de\,\{i^{-1}\}}_t\cap\{j\}\}}\big]\\[5pt]
\dis\quad=\int\!\chi_{\{0\}}P^\de_t\,(\di A)\,\int\!\chi_{\{0\}}\,(\di B)\,
|A\cap B|^{-1}1_{\{0\in A\cap B\}}\\[5pt]
\dis\quad=\lla\chi_{\{0\}}P^\de_t\captimes\chi_{\{0\}}\rra
=\lla\chi_{\{0\}}P^\de_s\captimes\chi_{\{0\}} P^{\dgg\,\de}_{t-s}\rra,
\ec
where we have used Lemma~\ref{L:infdual} in the last step.\qed

We will prove Theorem~\ref{T:difr} by taking the limit $t\to\infty$ in
(\ref{Russo2}). To justify the interchange of limit and differentiation, we
will use the following lemma.
\bl\pbf{(Interchange of limit and differentiation)}\label{L:diflim}
Let $I\sub\R$ be a compact interval and let $f_n,f,f'$ be continuous real
functions on $I$. Assume each $f_n$ is continuously differentiable, that
$f_n(x)\to f(x)$ and $\dif{x}f_n(x)\to f'(x)$ for each $x\in I$, and that
\be\label{supf}
\sup_{x\in I}\;\sup_n\,|\dif{x}f_n(x)|<\infty.
\ee
Then $f$ is continuously differentiable and $\dif{x}f(x)=f'(x)$ $(x\in I)$.
\el
\tbf{Proof} We write $I=[x_-,x_+]$ and observe that
\bc\label{intdif}
\dis f(x)&=&\dis\lim_{n\to\infty}f_n(x_-)
+\lim_{n\to\infty}\int_{x_-}^x\dif{y}f_n(y)\,\di y\\[5pt]
&=&\dis f(x_-)+\int_{x_-}^x\big(\lim_{n\to\infty}\dif{y}f_n(y)\big)\di y
=f(x_-)+\int_{x_-}^xf'(y)\,\di y,
\ec
where the interchange of limit and integration is justified by dominated
convergence, using (\ref{supf}). Differentiation of (\ref{intdif}) now yields
the statement since $f'$ is continuous.\qed

\noi 
\tbf{Proof of Theorem~\ref{T:difr}} Continuity of the map 
$(\de_{\rm c},\infty)\ni\de\mapsto\nucirc_\de$, and likewise for
$\nucirc^\dgg_\de$, in the sense of local convergence on $\Pfp$ has already
been proved in Proposition~\ref{P:eigpoint}. By Lemma~\ref{L:intsect}, this
implies local continuity of the map 
$(\de_{\rm c},\infty)\ni\de\mapsto\nucirc_\de\captimes\nucirc^\dgg_\de$.
Since local convergence on $\Pfp$ implies convergence of the
integral of the bounded functions $A\mapsto 1_{\{A=\{0\}\}}$ and
$A\mapsto |A|^{-1} 1_{\{0\in A\}}$ (which occurs in the definition of
$\lla\,\cdot\,\rra$), this implies continuity of
the right-hand side of (\ref{difr}).

Note that the right-hand side of (\ref{Russo}) is clearly bounded between zero
and one. Therefore, since 
\be
\frac{1}{t}\log\E\big[|\eta^{\de,\,\{0\}}_t|\big]
\asto{t}r(\La,a,\de)\qquad(\de\geq 0)
\ee
by the definition of the exponential growth rate in (\ref{rdef}), using
Lemma~\ref{L:diflim}, we see that (\ref{difr}) follows provided we show that
the right-hand side of (\ref{Russo2}) converges  for each
$\de\in(\de_{\rm c},\infty)$ to the right-hand side of (\ref{difr}) as $t \rightarrow \infty.$

We rewrite the right-hand side of (\ref{Russo2}) as
\be\label{uint}
\int_0^1\!\di u\,\frac{\dis 
e^{-rtu}\chi_{\{0\}}P^\de_{tu}
\captimes e^{-rt(1-u)}\chi_{\{0\}}P^{\dgg\,\de}_{t(1-u)}\,(\{0\})}
{\dis\lla e^{-rtu}\chi_{\{0\}}P^\de_{tu}\captimes
e^{-rt(1-u)}\chi_{\{0\}}P^{\dgg\,\de}_{t(1-u)}\rra}.
\ee
It is easy to see from the definition of $\lla\,\cdot\,\rra$ that the
integrand is bounded between zero and one (in fact, this is the probability in
(\ref{Russo})). By Theorem~\ref{T:eigcon}, for each $0<u<1$, the measures
$e^{-rtu}\chi_{\{0\}}P^\de_{tu}$ and
$e^{-rt(1-u)}\chi_{\{0\}}P^{\dgg\,\de}_{t(1-u)}$ converge locally on $\Pfp$ to
constant multiples of $\nucirc_\de$ and $\nucirc^\dgg_\de$, respectively. By
Lemma~\ref{L:intsect} and the fact that local convergence on $\Pfp$
implies convergence of the integral of the bounded functions
$A\mapsto 1_{\{A=\{0\}\}}$ and $A\mapsto |A|^{-1} 1_{\{0\in A\}}$, we see
that the integrand in (\ref{uint}) converges in a bounded pointwise way with
respect to $u$ to the right-hand side of (\ref{difr}).  Thus, the result
follows by Lebesgue's dominated convergence theorem.\qed

\section{Proof details}\label{S:detail}

In this section we supply the proof of all propositions and lemmas that have
not been proved yet. The organization is as follows. In Section~\ref{S:lf} we
prove some properties of locally finite measures and different forms of
convergence, concretely Proposition~\ref{P:locon} and Lemmas~\ref{L:intsect},
\ref{L:homfin}, \ref{L:cpsi} and \ref{L:timu}. In
Section~\ref{S:infstart} we consider contact processes started in infinite
initial `laws', proving Lemmas~\ref{L:infdual} and \ref{L:eiglim}. In
Section~\ref{S:expmom} we construct a metric on $\La$ with properties as in
Lemma~\ref{L:slowd} and prove Lemma~\ref{L:rga} on the exponential growth rate
associated with the functions $e_\ga$ defined in terms of such a metric. In
Section~\ref{S:covest} we do a covariance calculation leading to an estimate
of which Lemma~\ref{L:unmom} is a special case and use this to derive
Lemma~\ref{L:nutight}.

\subsection{Locally finite measures}\label{S:lf}

In this section, we prove Proposition~\ref{P:locon} as well as
Lemmas~\ref{L:intsect}, \ref{L:homfin} and \ref{L:cpsi}. Our first aim is
Proposition~\ref{P:locon}. We start with two preparatory lemmas. Recall the
definition of $\Pc_i$ from (\ref{Pidef}).

\bl\pbf{(Compact classes)}\label{L:Ccom}
If $\Ci\sub\Pc_+$ is compact, then there exists a finite $\De\sub\La$ such
that $\Ci\sub\bigcup_{i\in\De}\Pc_i$.
\el
\tbf{Proof} Choose $\De_n\up\La$ with $\De_n$ finite. If $\Ci\not\sub\bigcup_{i\in\De_n}\Pc_i$ for
each $n$, then we can find $A_n\in\Ci$ such that $A_n\cap\De_n=\emptyset$. It
follows that $A_n\to\emptyset\not\in\Ci$ (in the product topology),
hence $\Ci$ is not a closed subset of $\Pc$ and therefore not compact.\qed

\bl\pbf{(Vague and weak convergence)}\label{L:vagweak}
Let $\mu_n,\mu$ be locally finite measures on $\Pc_+$. Then the $\mu_n$
converge vaguely to $\mu$ if and only if for each $i\in\La$, the restricted
measures $\mu_n|_{\Pc_i}$ converge weakly to $\mu|_{\Pc_i}$ with respect to
the product topology.
\el
\tbf{Proof} Since $\Pc\beh\Pc_i$ is a closed subset
of $\Pc$, any continuous function $f:\Pc_i\to\R$ can be extended to a
continuous, compactly supported function on $\Pc_+$ by putting $f(A):=0$ for
$A\in\Pc_+\beh\Pc_i$. Therefore, if the $\mu_n$ converge vaguely to $\mu$, it
follows that the $\mu_n|_{\Pc_i}$ converge weakly to
$\mu|_{\Pc_i}$. Conversely, if for each $i\in\La$ the $\mu_n|_{\Pc_i}$
converge weakly to $\mu|_{\Pc_i}$, then for each $i,j\in\La$ one has
\be
\mu_n|_{\Pc_i\cap\Pc_j}\Rightarrow\mu|_{\Pc_i\cap\Pc_j},\quad
\mu_n|_{\Pc_i\beh\Pc_j}\Rightarrow\mu|_{\Pc_i\beh\Pc_j}\quad\mbox{and}\quad
\mu_n|_{\Pc_j\beh\Pc_i}\Rightarrow\mu|_{\Pc_j\beh\Pc_i},
\ee
where we have used that $\Pc_i\cap\Pc_j$, $\Pc_i\beh\Pc_j$ and
$\Pc_j\beh\Pc_i$ are compact sets. Continuing this process, we see by
induction that for each finite $\De\sub\La$, the restrictions
$\mu_n|_{\bigcup_{i\in\De}\Pc_i}$ converge weakly to
$\mu|_{\bigcup_{i\in\De}\Pc_i}$. By Lemma~\ref{L:Ccom}, if
$f:\Pc_+\to\R$ is a compactly supported continuous function, then $f$ is
supported on $\bigcup_{i\in\De}\Pc_i$ for some finite $\De\sub\La$. It follows
that $\int\mu_n(\di A)f(A)\to\int\mu(\di A)f(A)$, proving that the $\mu_n$
converge vaguely to $\mu$.\qed

\noi
\tbf{Proof of Proposition~\ref{P:locon}} The equivalence of (i) and (ii)
follows in a straightforward manner from Prohorov's theorem applied to the
countable space $\Pfi{i}$ with the discrete topology.


Since the discrete topology on $\Pfi{i}$ is stronger than the product
topology, weak convergence of the $\mu_n|_{\Pfi{i}}$ with respect to the
discrete topology implies weak convergence with respect to the product
topology. By Lemma~\ref{L:vagweak}, this shows that local convergence on
$\Pfp$ implies vague convergence on $\Pc_+$ and hence (i) implies also
(iii).

To prove (iii)$\Rightarrow$(i), note that by local tightness, for each
$i\in\La$ the measures $\mu_n|_{\Pfi{i}}$ are relatively compact in the
topology of weak convergence with respect to the discrete topology. Let
$\mu^i_\ast$ be a subsequential limit. Since weak convergence with respect to
the discrete topology implies weak convergence with respect to the product
topology, by Lemma~\ref{L:vagweak}, we conclude that
$\mu^i_\ast=\mu|_{\Pfi{i}}$. Since this is true for each cluster point, we
conclude that the $\mu_n|_{\Pfi{i}}$ converge weakly to $\mu|_{\Pfi{i}}$ with
respect to the discrete topology.

The implication (i)$\volgt$(iv) follows from what we have already proved. To
prove the reverse implication, it suffices to show local tightness. Since for
each $i\in\La$, the finite measures $\mu_n|_{\Pfi{i}}$ converge pointwise to
$\mu|_{\Pfi{i}}$, it suffices to show that their total mass satisfies
\be
\limsup_{n\to\infty}\mu_n(\{A:i\in A\})\leq\mu(\{A:i\in A\}).
\ee
By vague convergence (see Lemma~\ref{L:locfin}), the limit superior is
actually a limit and equals the right-hand side.\qed

\noi
\tbf{Proof of Lemma~\ref{L:intsect}} The local finiteness of $\mu\captimes\nu$
follows from Lemma~\ref{L:locfin} and the fact that
\bc
\dis\int\mu\captimes\nu\,(\di C)1_{\{i\in C\}}
&=&\dis\int\mu(\di A)\int\nu(\di B)1_{\{i\in A\cap B\}}\\[5pt]
&=&\dis\Big(\int\mu(\di A)1_{\{i\in A\}}\Big)
\Big(\int\nu(\di B)1_{\{i\in B\}}\Big)<\infty\quad(i\in\La).
\ec
To see that $\mu_n\captimes\nu_n$ converges vaguely to $\mu\captimes\nu$ if
$\mu_n,\nu_n$ converge vaguely to $\mu,\nu$, respectively, by
Lemma~\ref{L:locfin}, it suffices to check that
\be\label{psico}
\int\mu_n\captimes\nu_n(\di C)1_{\{C\cap D\neq\emptyset\}}
\asto{n}\int\mu\captimes\nu(\di C)1_{\{C\cap D\neq\emptyset\}}
\qquad(D\in\Pfp).
\ee
Since
\be
1_{\{C\cap D\neq\emptyset\}}
=1-\prod_{i\in D}1_{\{i\not\in C\}}=1-\prod_{i\in D}(1-1_{\{i\in C\}})
=\sum\subb{D'\sub D}{D'\neq\emptyset}(-1)^{|D'|+1}\prod_{i\in D'}1_{\{i\in C\}},
\ee
and since $\prod_{i\in D'}1_{\{i\in C\}} = 1_{\{D' \subset C\} }$ formula (\ref{psico}) is equivalent to
\be
\int\mu_n\captimes\nu_n(\di C)1_{\{D\sub C\}}
\asto{n}\int\mu\captimes\nu(\di C)1_{\{D\sub C\}}
\qquad(D\in\Pfp).
\ee
Now
\bc
\dis\int\mu_n\captimes\nu_n(\di C)1_{\{D\sub C\}}
&=&\dis\int\mu_n(\di A)\int\nu_n(\di B)
1_{\{D\sub(A\cap B)\}}\\[5pt]
&=&\dis\Big(\int\mu_n(\di A)1_{\{D\sub A\}}\Big)
\Big(\int\nu_n(\di B)1_{\{D\sub B\}}\Big),
\ec
which, by our assumptions that $\mu_n\Rightarrow\mu$ and $\nu_n\Rightarrow\nu$,
converges to the analogue formula with $\mu_n,\nu_n$ replaced by $\mu,\nu$.

To see that the vague convergence of $\mu_n\captimes\nu_n$ can be strengthened
to local convergence on $\Pfp$ if either $\mu_n$ or $\nu_n$ converges locally
on $\Pfp$,  it suffices by Proposition~\ref{P:locon}~(iii)$\Rightarrow$(i) to show that the
local tightness of either $\mu_n$ or $\nu_n$ implies local tightness of 
$\mu_n\captimes\nu_n$. By symmetry, it suffices to consider the case when the
$\mu_n$ are locally tight. Since vague convergence of the $\nu_n$ implies
convergence of $\int\nu_n(\di A)1_{\{i\in A\}}$ for each $i\in\La$, the
statement now follows from the following lemma, that we formulate separately
since it is of some interest on its own.\qed

\bl\pbf{(Local tightness of intersection measure)}\label{L:captight}
Let $\mu_n,\nu_n$ $(n\geq 1)$ be locally finite measures on $\Pc_+$. Assume
that the $\mu_n$ $(n\geq 1)$ are concentrated on $\Pfp$ and that they are locally tight. Assume
that the $\nu_n$ satisfy $\sup_{n\geq 1}\int\nu_n(\di A)1_{\{i\in A\}}<\infty$
for all $i\in\La$. Then the intersection measures $\mu_n\captimes\nu_n$
$(n\geq 1)$ are concentrated on $\Pfp$ and locally tight.
\el
\tbf{Proof} Since $\mu_n\captimes\nu_n$ is concentrated on sets of the form
$A\cap B$ with $A\in\Pfp$, it is clear that $\mu_n\captimes\nu_n$ is
concentrated on $\Pfp$ for each $n\geq 1$. Fix $i\in\La$ and $\eps>0$,
and set $K:=\sup_{n\geq 1}\int\nu_n(\di A)1_{\{i\in A\}}$. By the local tightness of the $\mu_n$, there exists a finite
$\Di\sub\Pfi{i}$ such that $\sup_n\mu_n(\Pfi{i}\beh\Di)\leq\eps/K$.
The same obviously holds for the larger finite set $\Di':=\Pc(D)=\{A:A\sub
D\}$, where $D:=\bigcup\{A:A\in\Di\}$. Now
\be\ba{l}
\dis\sup_{n\geq 1}\mu_n\captimes\nu_n(\Pfi{i}\beh\Di')
=\sup_{n\geq 1}\int\!\mu_n(\di A)\,\int\!\nu_n(\di B)\,
1_{\{i\in A\cap B\}}1_{\{A\cap B\not\sub D\}}\\[5pt]
\dis\quad\leq\sup_{n\geq 1}\int\!\mu_n(\di A)\,1_{\{i\in A\}}1_{\{A\not\sub D\}}
\int\!\nu_n(\di B)\,1_{\{i\in B\}}\leq\eps.
\ec
Since $i\in\La$ and $\eps>0$ are arbitrary, the claim follows.\qed

\noi
\tbf{Proof of Lemmas~\ref{L:homfin} and \ref{L:timu}}
Formula (\ref{muDe}) obviously defines a nonzero, homogeneous measure on
$\Pfp$. Since
\be
\mu(\{A:0\in A\})=c\sum_i\P[0\in i\De\big]=c\sum_i\P[i^{-1}\in\De\big]
=cE\big[|\De|\big],
\ee
it follows from Lemma~\ref{L:locfin} that $\mu$ is locally finite
if and only if $\E\big[|\De|\big]<\infty$. If $\mu$ is given by (\ref{muDe}),
then
\be
\lla\mu\rra
=c\sum_{i\in\La}\E\big[|i\De|^{-1}1_{\{0\in i\De\}}\big]
=c\E\big[|\De|^{-1}\big(\sum_{i\in\La}1_{\{i^{-1}\in\De\}}\big)\big]=c.
\ee
To see that every nonzero, homogeneous measure $\mu$ on $\Pfp$ with
$\lla\mu\rra<\infty$ can be written in the form (\ref{muDe}), define a
probability law $\rho$ on $\Pfi{0}$ by
\be
\rho(\{A\}):=\lla\mu\rra^{-1}\mu(\{A\})|A|^{-1}1_{\{0\in A\}}.
\ee
Let $\De$ be a random variable with law $\rho$. We claim that $\mu$ is given
by (\ref{muDe}) with $c=\lla\mu\rra$. To check this, we calculate, for
$A\in\Pfp$:
\be\ba{l}
\dis \lla\mu\rra\sum_{i\in\La}\P\big[i\De=A\big]
=\lla\mu\rra\sum_{i\in\La}\P\big[\De=i^{-1}A\big]
=\lla\mu\rra\sum_{i\in\La}\rho(\{i^{-1}A\})\\[5pt]
\dis\quad=\sum_{i\in\La}\mu(\{ i^{-1} A \})|i^{-1} A|^{-1}1_{\{0\in i^{-1}A\}}
=\mu(\{A\})|A|^{-1}\sum_{i\in\La}1_{\{i\in A\}}=\mu(\{A\}),
\ec
where we have used the homogeneity of $\mu$. This completes the proof of
Lemma~\ref{L:homfin}, except for the statement that the law of $\ti\De$ is
uniquely determined by $\mu$, which will follow by setting $t=0$ in formula
(\ref{mutimu}) of Lemma~\ref{L:timu}, which we prove next.

Indeed, letting $\eta_0$ be a $\Pfp$-valued random variable such that
\be
\mu=\lla\mu\rra\sum_{i\in\La}\P[i\eta_0\in\cdot\,],
\ee
and letting $(\eta_t)_{t\geq 0}$ be a $(\La,a,\de)$-contact process started in
$\eta_0$, we have for any $B\in\Pfp$ and $t\geq 0$ that
\be\ba{l}\label{mup}
\dis\mu P_t(\{B\})=\!\sum_{A\in\Pfp}\!\mu(\{A\})P_t(A,B)
=\lla\mu\rra\!\sum_{A\in\Pfp}\!\sum_{i\in\La}\P[i\eta_0=A]P_t(A,B)\\[5pt]
\dis\quad=\lla\mu\rra\!\sum_{A\in\Pfp}\!\sum_{i\in\La}\P[i\eta_t=B]
=\lla\mu\rra m(B)\P[\ti\eta_t=\ti B],
\ec
where $m(B)$ is defined as in (\ref{mA}). Letting $\ti\mu$ denote the law of
$\ti\eta_0$ and $\ti P_t$ the transition probabilities of the
$(\La,a,\de)$-contact process modulo shifts, we arrive from (\ref{mup} at
(\ref{mutimu}).

To complete also the proof of Lemma~\ref{L:timu}, we still need to prove
(\ref{timuf}). Let $\acPfp$ be a set that contains exactly one representative
from each equivalence class $\ti A\in\tiPfp$. Then (\ref{timuf}) follows,
using (\ref{mutimu}), by writing
\be\ba{l}
\dis\lla\mu\rra\!\sum_{\ti A\in\tiPfp}\!\ti\mu(\ti A)\ti f(\ti A)
=\!\sum_{A\in\acPfp}\!m(A)^{-1}\mu(\{A\})f(A)|A|^{-1}\sum_{i\in\La}1_{\{i\in A\}}\\[5pt]
\dis\quad=\!\sum_{A\in\acPfp}\!m(A)^{-1}\sum_{i\in\La}(f\mu)(\{i^{-1}A\})
|i^{-1}A|^{-1}1_{\{0\in i^{-1}A\}}=\lla f\mu\rra.
\ec
\qed



\noi
We finish the section on locally finite measures by supplying the still
outstanding:\med

\noi
\tbf{Proof of Lemma~\ref{L:cpsi}} We will apply the mass transport principle,
compare the proof of Lemma~\ref{L:infdual} below. Let $\mu,\nu$ be
homogeneous, locally finite measures on $\Pc_+$ and assume that $\mu$ is
concentrated on $\Pfp$. For $A\in\Pf$ and $B\in\Pc$ such that $A\cap
B\neq\emptyset$, let us define a probability distribution $M_{A,B}$ on
$\La\times\La$ by
\be
M_{A,B}(i,j)
:=|A|^{-1}1_{\{i\in A\}}|A\cap B|^{-1}1_{\{j\in A\cap B\}},
\ee
and let $f:\La\times\La\to[0,\infty]$ be defined by
\be
f(i,j):=\int\!\mu(\di A)\,\int\!\nu(\di B)\,
1_{\{A\cap B\neq\emptyset\}}M_{A,B}(i,j).
\ee
Since $\mu$ and $\nu$ are homogeneous, we observe that $f(ki,kj)=f(i,j)$
$(i,j,k\in\La)$. Moreover,
\be
\sum_jf(0,j)=\int\!\mu(\di A)\,\int\!\nu(\di B)\,1_{\{A\cap B\neq\emptyset\}}
\frac{1}{|A|}1_{\{0\in A\}}
=\int\!\mu(\di A)\,h_\nu(A)\frac{1}{|A|}1_{\{0\in A\}}=\lla h_\nu \mu\rra,
\ee
while
\be
\sum_if(i,0)=\int\!\mu(\di A)\,\int\!\nu(\di B)\,1_{\{A\cap B\neq\emptyset\}}
\frac{1}{|A\cap B|}1_{\{0\in A\cap B\}}
=\lla\mu\captimes\nu\rra.
\ee
Formula (\ref{cpsi}) now follows from the fact that
$\sum_if(i,0)=\sum_if(0,i^{-1})=\sum_jf(0,j)$.
Note that this holds regardless of whether $h_\nu\mu$ is locally finite or
not. If $\int\mu(\di A)|A|1_{\{0\in A\}}<\infty$, then by the
shift-invariance and subadditivity of $h_\nu$, we see that $h_\nu(A)\leq
h_\nu(\{0\})|A|$ and hence $\int\mu(\di A)h_\nu(A)1_{\{0\in A\}}<\infty$,
proving that $h_\nu\mu$ is locally finite.\qed

\subsection{Infinite starting measures}\label{S:infstart}

In this section we prove Lemma~\ref{L:infdual} on contact process duality
for homogeneous, infinite starting measures. We also give
the proof of Lemma~\ref{L:eiglim}, which is concerned with relative
compactness and cluster points of eigenmeasures for $(\La,a,\de)$-contact
processes with varying~$\de$.\med

\noi 
\tbf{Proof of Lemma~\ref{L:infdual}} Fix $t\geq 0$ and for $A,B\in\Pc_+$,
consider the events
\be
\Ei_{A,B}:=\{|\eta^{A,0}_t\cap B|<\infty\}
\quad\mbox{and}\quad
\Ei'_{A,B}:=\{|A\cap\eta^{\dgg\,B,t}_t|<\infty\}.
\ee
We observe that $\mu P_t\captimes\nu$ (resp.\ $\mu\captimes\nu P^\dgg_t$) is
concentrated on $\Pfp$ if and only if $\P(\Ei_{A,B})=1$
(resp.\ $\P(\Ei'_{A,B})=1$) for a.e.\ $A$ w.r.t.\ $\mu$ and a.e.\ $B$
w.r.t.\ $\nu$. Set $\De_0:=A\cap\eta^{\dgg\,B,t}_t$ and
$\De_t:=\eta^{A,0}_t\cap B$. Since $\eta^{\De_0,0}_t\supset\De_t$ and
$\eta^{\dgg\,\De_t,t}_t\supset\De_0$, we see that the events $\Ei_{A,B}$ and
$\Ei'_{A,B}$ are a.s.\ equal, and hence $\mu P_t\captimes\nu$ is
concentrated on $\Pfp$ if and only if $\mu\captimes\nu P^\dgg_t$ is.

We will now prove (\ref{infdual}) by applying the ``mass transport
principle''. For a given graphical representation $\om$ and
sets $A,B\in\Pc_+$ such that the events $\Ei_{A,B}$ and $\Ei'_{A,B}$ hold, we
define a probability distribution $M_{A,B,\om}$ on $\La\times\La$ by
\be
M_{A,B,\om}(i,j):=|\De_0|^{-1}1_{\{i\in\De_0\}}|\De_t|^{-1}1_{\{j\in\De_t\}}.
\ee
We define a function $f:\La\times\La\to[0,\infty]$ by
\be
f(i,j):=\int\!\mu(\di A)\int\!\nu(\di B)\,\int\!\P(\di\om)\,
1_{\Ei_{A,B}}(\om)\,M_{A,B,\om}(i,j).
\ee
Obviously, $f(ki,kj)=f(i,j)$ $(i,j,k\in\La)$ due to the homogeneity of $\mu$
and $\nu.$ Moreover,
\be\ba{l}
\dis\sum_if(i,0)=\int\!\mu(\di A)\int\!\nu(\di B)\,
\E\big[|\eta^{A,0}_t\cap B|^{-1}1_{\{0\in\eta^{A,0}_t\cap B\}}\big]\\[5pt]
\dis\quad=\int\mu P_t(\di A')\int\nu(\di B)
|A'\cap B|^{-1}1_{\{0\in A'\cap B\}}
=\lla\mu P_t\captimes\nu\rra.
\ec
The same argument shows that $\sum_jf(0,j)=\lla\mu\captimes\nu P^\dgg_t\rra$
and hence
\be
\lla\mu P_t\captimes\nu\rra=\sum_if(i,0)=\sum_if(0,i^{-1})
=\lla\mu\captimes\nu P^\dgg_t\rra,
\ee
where the middle step is a simple example of what is more generally known as
the mass transport principle, see \cite{Hag11}.\qed

\noi
\tbf{Proof of Lemma~\ref{L:eiglim}} By the homogeneity and normalization of
the $\nu_n$, one has
\be\label{Babs}
\int\nu_n(\di A)1_{\{A\cap B\neq\emptyset\}}
\leq\sum_{i\in B}\int\nu_n(\di A)1_{\{i\in A\}}=|B|.
\ee
Since this estimate is uniform in $n$, applying \cite[Lemma~3.2]{Swa09} we
find that the $(\nu_n)_{n\geq 0}$ are relatively compact in the topology of
vague convergence. By going to a subsequence if necessary, we may assume that
the $\nu_n$ converge vaguely to a limit $\nu$. Since the $\nu_n$ are
eigenmeasures, denoting the $(\La,a,\de_n)$-contact process started in $A$ by
$(\eta^{\de_n,A}_t)_{t\geq 0}$, we have
\be\label{nu_neig}
\int\nu_n(\di A)\P[\eta^{\de_n,A}_t\in\cdot\,]\big|_{\Pc_+}
=e^{\la_nt}\nu_n\qquad(t\geq 0).
\ee
Since $\la_n\to\la$, the right-hand side of this equation converges vaguely to
$e^{\la t}\nu$. To prove vague convergence of the left-hand side, by
Lemma~\ref{L:locfin}, it suffices to prove that for $B \in \Pf$,
\be\label{nu_neigleftconv}
\int\nu_n(\di A)\P[\eta^{\de_n,A}_t\cap B\neq\emptyset]\rightarrow
\int\nu(\di A)\P[\eta^{\de,A}_t\cap B \neq\emptyset].
\ee
We estimate
\begin{eqnarray}
\nonumber
&&\Big|\int\nu_n(\di A)\P[\eta^{\de_n,A}_t\cap B\neq\emptyset]
-\int\nu(\di A)\P[\eta^{\de,A}_t\cap B \neq\emptyset]\Big|\\
\label{nu_nconv1}
&&\qquad\leq\int\nu_n(\di A)\Big|\P[\eta^{\de_n,A}_t\cap B\neq\emptyset]
-\P[\eta^{\de,A}_t\cap B \neq\emptyset]\Big|\\
\label{nu_nconv2}
&&\dis\qquad\phantom{\leq}+\Big|\int\nu_n(\di A)
\P[\eta^{\de,A}_t\cap B\neq\emptyset]
-\int\nu(\di A)\P[\eta^{\de,A}_t\cap B \neq\emptyset]\Big|.
\end{eqnarray}
The term in (\ref{nu_nconv2}) tends to zero as $n\to\infty$ by Lemmas
\ref{L:locfin} and \ref{L:infstart}. By duality, we can rewrite the term in
(\ref{nu_nconv1}) as
\be\label{du_nconv1}
\int\nu_n(\di A)\Big|\P[A\cap\eta^{\dgg\,\de_n,B}_t\neq\emptyset]
-\P[A\cap\eta^{\dgg\,\de,B}_t\neq\emptyset]\Big|.
\ee
We couple the graphical representations for processes with different recovery
rates in the natural way, by constructing a Poisson point process $\Om^{\rm
r}$ on $\La\times\R_+\times\R_+$ with intensity one, and letting $\poi^{\rm
r}_\de:=\{(i,t):\exists 0\leq r\leq\de\mbox{ s.t.\ }(i,t,r)\in\Om^{\rm r}\}$
be the set of recovery symbols for the process with recovery rate $\de$.
Then, letting $\eta^{\dgg\,0,B}_t$ denote the process with zero recovery rate,
the quantity in (\ref{du_nconv1}) can be estimated from above by
\be\ba{l}\label{etanul}
\dis\int\nu_n(\di A)\P\big[A\cap\eta^{\dgg\,0,B}_t\neq\emptyset,
\ \eta^{\dgg\,\de_n,B}_t\neq\eta^{\dgg\,\de,B}_t\big]\\[5pt]
\dis\quad=\int\P\big[\eta^{\dgg\,0,B}_t\in\di C,
\ \eta^{\dgg\,\de_n,B}_t\neq\eta^{\dgg\,\de,B}_t\big]
\int\nu_n(\di A)1_{\{A\cap C\neq\emptyset\}}\\[5pt]
\dis\quad\leq\int\P\big[\eta^{\dgg\,0,B}_t\in\di C,
\ \eta^{\dgg\,\de_n,B}_t\neq\eta^{\dgg\,\de,B}_t\big]\,|C|
=\E\big[|\eta^{\dgg\,0,B}_t|
1_{\{\eta^{\dgg\,\de_n,B}_t\neq\eta^{\dgg\,\de,B}_t\}}\big],
\ec
where we have used (\ref{Babs}). Since the right-hand side of (\ref{etanul})
tends to zero by dominated convergence, this proves the lemma.\qed

\subsection{Exponential moments}\label{S:expmom}

Recall the function $e_{\gamma}(A) = \sum_{i \in A} e^{\gamma d(0,i)}$ from
(\ref{ega}), which measures how `spread out' a set $A \in \Pf$ is in terms of
exponential weights and a suitably slowly growing metric $d$ as in
(\ref{dprop}). In this section, we provide the proof of Lemma~\ref{L:slowd},
showing that such a metric exists. We then give the proof of
Lemma~\ref{L:rga}, which states that the expectation of the function $e_{\ga}$
of a contact process has a well defined exponential growth rate, with certain
bounds.\med

\noi
\tbf{Proof of Lemma~\ref{L:slowd}} We can find finite
$\{0\}=\De_1\sub\De_2\sub\cdots$ such that
$\sum_{i\in\La\beh\De_n}a(0,i)\leq|a|e^{-(n-1)}$. Making the sets $\De_n$ for
$n\geq 2$ larger if necessary, we can moreover choose these sets such that
they are symmetric, i.e., $\{i^{-1}:i\in\De_n\}=\De_n$ and such that
$\De_\infty:=\bigcup_{n\geq 1}\De_n$ generates $\La$. (In particular, we can
always choose $\De_\infty=\La$, but for nearest-neighbor processes on graphs
this leads to a somewhat unnatural metric $d$, which is why we only assume
here that $\De_\infty$ generates $\La$.) We set $\De_0:=\emptyset$ and define
\be
\phi(i):=\left\{\ba{ll}
n\qquad&(i\in\De_n\beh\De_{n-1},\ n\geq 1)\\[5pt]
\infty\qquad&(i\in\La\beh\De_\infty).
\ea\right.\ee
Since $a(0,i)=0$ for $i\not\in\De_\infty$, whe have that 
\be
\sum_{i\in\La}a(0,i)\phi(i)^\ga
=\sum_{n\geq 1}n^\ga\sum_{i\in\De_n\beh\De_{n-1}}a(0,i)
\leq|a|\sum_{n\geq 1}n^\ga e^{-(n-2)}<\infty
\ee
for each $0\leq\ga<\infty$. Set
\be
d'(i,j)=d'(0,i^{-1}j):=\log(\phi(i^{-1}j))\qquad(i,j\in\La).
\ee
Then $d'$ satisfies properties (\ref{dprop})~(i)--(iii), $d'(i,j)=0$ if and
only if $i=j$, and $d'(i,j)=d'(j,i)$ (by the symmetry of the sets $\De_n$).
Since $d'$ need not yet be a metric, we define
\be
d(i,j):=\inf\big\{\sum_{k=1}^nd'(i_{k-1},i_k):
n\geq 1,\ i_0,\ldots,i_n\in\La,\ i_0=i,\ i_n=j\big\},
\ee
i.e., $d(i,j)$ is a graph-style distance between $i$ and $j$, defined as the
shortest path from $i$ to $j$ where an edge from $i_{k-1}$ to $i_k$ has length
$d'(i_{k-1},i_k)$. Note that $d(i,j)<\infty$ for each $i,j\in\La$ since
$\De_\infty$ generates $\La$ and $d(i,j)>0$ for each $i\neq j$ since
$d'(i,j)\geq\log(2)$ for each $i\neq j$. It is now straightforward to check
that $d$ is a metric on $\La$ and that $d(i,j)=d(ki,kj)$ for all
$i,j,k\in\La$. Since $d(i,j)\leq d'(i,j)$, the metric $d$ also enjoys property
(\ref{dprop})~(iii). Property (\ref{dprop})~(ii), finally, follows
from the fact that
\be
\{i\in\La:d(0,i)\leq M\}
\sub\{j_1\cdots j_n:1\leq n\leq M/\log(2),
\ d'(0,j_k)\leq M\ \forall k=1,\ldots,n\},
\ee
where we use that $d'(i,j)\geq\log(2)$ for all $i\neq j$, and we observe
that if $d(0,i)\leq M$ $(i\neq 0)$, then there must be some
$n\geq 1$ and $0=i_0,\ldots,i_n=i$ with $\sum_{k=1}^nd'(i_{k-1},i_k)\leq
M$. Setting $j_k:=i_{k-1}^{-1}i_k$ we see that $i$ must be of the form
$i=j_1\cdots j_n$ with $\sum_{k=1}^nd'(0,j_k)\leq M$.\qed

As a preparation for Lemma~\ref{L:rga}, we need one more result.

\bl\pbf{(Existence of exponential moments)}\label{L:expmom}
Let $(\eta^A_t)_{t\geq 0}$ be a $(\La,a,\de)$-contact process started in a
finite initial state $\eta^A_0=A\in\Pf$ and let $d$ be a metric on $\La$ as in
Lemma~\ref{L:slowd}. Then
\be\label{Kgdef}
\E\big[e_\ga(\eta^A_t)\big]\leq e^{K_\ga t}e_\ga(A)\quad(t\geq 0)
\quad\mbox{where}\quad
K_\ga:=\sum_{i\in\La}a(0,i)e^{\ga d(0,i)}.
\ee
\el
\tbf{Proof} For $\ga=0$ this follows from \cite[Prop.~2.1]{Swa09}. To prove
the statement for $\ga>0$, let $G$ be the generator of the
$(\La,a,\de)$-contact process as defined in (\ref{Gdef}). Then
\bc\label{Kge}
\dis Ge_\ga(A)&=&\dis\sum_{i\in A}\sum_{j\not\in A}a(i,j)e^{\ga d(0,j)}
-\de\sum_{i\in A}e^{-\ga d(0,i)}\\[5pt]
&\leq&\dis\sum_{i\in A}\sum_{j\in\La}a(i,j)e^{\ga(d(0,i)+d(i,j))}
=K_\ga e_\ga(A),
\ec
where we have used that
$\sum_{j\in\La}a(i,j)e^{\ga d(i,j)}
=\sum_{j\in\La}a(0,i^{-1}j)e^{\ga d(0,i^{-1}j)}=K_\ga$
$(i\in\La)$.

Set $\tau_N:=\inf\{t\geq 0:e_\ga(\eta^A_t)\geq N\}$. Since the stopped process
is a Markov process with finite state space, it follows by standard arguments
from (\ref{Kge}) that
\be\label{stopest}
\E\big[e_\ga(\eta^A_{t\wedge\tau_N})\big]\leq e^{K_\ga t}e_\ga(A)
\quad(t\geq 0,\ N\geq 1),
\ee
which in turn implies that $\P[e_\ga(\eta^A_{t\wedge\tau_N})\geq N]\to 0$ as
$N\to\infty$ and hence $\tau_N\to\infty$ a.s. Therefore, letting $N\to\infty$
in (\ref{stopest}), we arrive at (\ref{Kgdef}).\qed

\noi
\tbf{Proof of Lemma~\ref{L:rga}} Note that $r_0(\La,a,\de)=r(\La,a,\de)$ is
the exponential growth rate from (\ref{rdef}). The statement for $\ga=0$ has
been proved in \cite[Lemma~1.1 and formula (3.5)]{Swa09}. To prove the general
statement, set $\pi^\ga_t:=\E\big[e_\ga(\eta^{\{0\}}_t)\big]$. Formula
(\ref{rga}) will follow from standard facts \cite[Thm~B.22]{Lig99} if we show
that $t\mapsto\log\pi^\ga_t$ is subadditive. Recalling the graphical
representation of the $(\La,a,\de)$-contact process, we observe that indeed
\be\ba{l}
\dis\pi^\ga_{s+t}=\sum_i\P[(0,0)\leadsto(i,s+t)]e^{\ga d(0,i)}\\[5pt]
\dis\quad\leq
\sum_{ij}\P[(0,0)\leadsto(j,s)\leadsto(i,s+t)]e^{\ga(d(0,j)+d(j,i))}
=\pi^\ga_s\pi^\ga_t,
\ec
which implies the subadditivity of $t\mapsto\log\pi^\ga_t$ and hence formula
(\ref{rga}). Since $e_\ga(A)\leq e_{\ga'}(A)$ for all $\ga\leq\ga'$, it is
clear that $\ga\mapsto r_\ga$ is nondecreasing. The fact that $-\de\leq r_0$
has been proved in \cite[Lemma~1.1]{Swa09} while the estimate $r_\ga\leq
K_\ga$ is immediate from Lemma~\ref{L:expmom}.

To prove that the function $\half\ni\ga\mapsto r_\ga$ defined in
Lemma~\ref{L:rga} is right-continuous, we observe that it follows from
(\ref{rga}) that for any $t_n\up\infty$,
\be
r_\ga=\lim_{n\to\infty}\inf_{1\leq k\leq n}
\frac{1}{t_k}\log\E\big[e_\ga(\eta^{\{0\}}_{t_k})\big].
\ee
By dominated convergence and the finiteness of exponential moments
(Lemma~\ref{L:expmom}) we have that for each fixed $t>0$, the function
$\ga\mapsto\frac{1}{t}\log\E[e_\ga(\eta^{\{0\}}_t)]$ is continuous. Therefore,
being the decreasing limit of continuous functions, $\ga\mapsto r_\ga$ must be
upper semi-continuous. Since $\ga\mapsto r_\ga$ is nondecreasing, this is
equivalent to continuity from the right.\qed

\subsection{Covariance estimates}\label{S:covest}

The next lemma gives a uniform estimate on expectations of the functions
$e_\ga(A)$ defined in (\ref{ega}) under the measures $1_{\{0\in \cdot\}}
\frac{1}{\hat\pi_\la}\hat\mu_\la$. Lemma~\ref{L:unmom} and
Lemma~\ref{L:nutight}, which were stated and used in
Sections~\ref{S:Existence} and \ref{S:cont} respectively, follow as
corollaries to this lemma. Their proofs are given at the end of this section.

Although this is not exactly how the proof goes, the following heuristic is
perhaps useful for understanding the main strategy. Since Campbell measures
change second moments into first moments, what we need to control are
expectations of the form $\E[e_\ga(\eta^{\{0\}}_t)^2]$, which leads us to
consider events of the form
\be
(0,0)\leadsto(i,t)\quad\mbox{and}\quad(0,0)\leadsto(j,t).
\ee
Since in the subcritical regime, long connections are unlikely, the
largest contribution to the probability of such an event comes from events of
the form
\be\label{vias}
(0,0)\leadsto(k,s)\left\{\ba{l}\leadsto(i,t)\\ \leadsto(j,t).\ea\right. 
\ee
where $s\in[0,t]$ is close to $t$ and $k\in\La$. Indeed, if the exponential
growth rate $r=r(\La,a,\de)$ is negative, then the probability of an event of
the form (\ref{vias}) is of the order $e^{rs}(e^{r(t-s)})^2$, which much
smaller than the probability that $(0,0)\leadsto(i,t)$, unless $t-s$ is of
order one. In view of this, if we find an infection at some late time $t$,
then all other infected sites are likely to be close to it. Although this
reasoning is only heuristic, it turns out that the covariance formula
(\ref{varfor}) below provides a convenient way of making such arguments
precise.

\bl\pbf{(Uniform exponential moment bound)}\label{L:unexmom}
Let $\hat\mu_\la$ and $\hat\pi_\la$ be defined as in
(\ref{Lap})--(\ref{hatpi}) and for $\ga\geq 0$, let $e_\ga$ be the function
defined in (\ref{ega}) in terms of a metric $d$ satisfying (\ref{dprop}).
Then, for any $(\La,a,\de)$-contact process with exponential growth rate
$r=r(\La,a,\de)$,
\be\label{unexmom}
\limsup_{\la\down r}\frac{1}{\hat\pi_\la}
\int\hat\mu_\la(\di A)1_{\{0\in A\}}e_\ga(A)
\leq(|a|+\de)\int_0^\infty e^{-rt}\di t\,
\E\big[e_\ga(\eta^{\{0\}}_t)\big]^2.
\ee
\el
We note that although the bound
in (\ref{unexmom}) holds regardless of the values of $\ga$ and
$r=r(\La,a,\de)$, the right-hand side will usually be infinite, unless $r<0$
and $\ga$ is small enough (see the proofs of Lemma~\ref{L:finex} and
Proposition~\ref{P:eigpoint}).\med

\noindent  
\tbf{Proof} Fix $\ga\geq 0$ and, to ease notation, set $\psi_\ga(i,j):=e^{\ga
d(i,j)}$ $(i,j,k\in\La)$. We observe that
\be\ba{l}\label{mud}
\dis\int\hat\mu_\la(\di A)1_{\{0\in A\}}e_\ga(A)
=\int_0^\infty e^{-\la t}\di t\,\sum_{i,j}
\E\big[1_{\{0\in\eta^{\{i\}}_t\}}1_{\{j\in\eta^{\{i\}}_t\}}
\psi_\ga(0,j)\big]\\[5pt]
\dis\quad=\int_0^\infty e^{-\la t}\di t\,\sum_{i,j}
\E\big[1_{\{i^{-1}\in\eta^{\{0\}}_t\}}1_{\{i^{-1}j\in\eta^{\{0\}}_t\}}
\psi_\ga(i^{-1},i^{-1}j)\big]\\[5pt]
\dis\quad=\int_0^\infty e^{-\la t}\di t\,\sum_{i,j}\psi_\ga(i,j)
\P\big[i\in\eta^{\{0\}}_t,\ j\in\eta^{\{0\}}_t\big].
\ec
Set $f_i(A):=1_{\{i\in A\}}$. Then
\be\label{fifj}
\P\big[i\in\eta^{\{0\}}_t,\ j\in\eta^{\{0\}}_t\big]
=\E\big[f_i(\eta^{\{0\}}_t)\big]\,\E\big[f_j(\eta^{\{0\}}_t)\big]
+\cov\big(f_i(\eta^{\{0\}}_t),f_j(\eta^{\{0\}}_t)\big).
\ee
By a standard covariance formula (see \cite[Prop.~2.2]{Swa09}), for any
functions $f,g$ of polynomial growth (as in (\ref{Gdom}) below), one has
\be\label{varfor}
\cov\big(f(\eta^{\{0\}}_t),g(\eta^{\{0\}}_t)\big)
=2\int_0^t\E\big[\Ga(P_sf,P_sg)(\eta^{\{0\}}_{t-s})\big]\di s
\qquad(t\geq 0),
\ee
where $(P_t)_{t\geq 0}$ denotes the semigroup of the $(\La,a,\de)$-contact
process and $\Ga(f,g)=\ffrac{1}{2}(G(fg)-fGg-gGf)$, with $G$ as in
(\ref{Gdef}). A little calculation (see \cite[formula (4.6)]{Swa09}) shows that
\bc\label{2Gamma}
\dis 2\Ga(P_sf,P_sg)(A)
&=&\dis\sum_{k\in A}\sum_{l\not\in A}a(k,l)
\big(P_sf(A\cup\{l\})-P_sf(A)\big)(P_sg(A\cup\{l\})-P_sg(A)\big)\\[5pt]
&&\dis+\de\sum_{k\in A}
\big(P_sf(A\beh\{k\})-P_sf(A)\big)\big(P_sg(A\beh\{k\})-P_sg(A)\big).
\ec
Applying (\ref{2Gamma}) to the functions $f=f_i$, $g=f_j$, using the fact
that, by the graphical representation,
\be
\big|P_sf_i(A\cup\{l\})-P_sf_i(A)\big|
=\big|\P\big[i\in\eta^{A\cup\{l\}}_s]-\P\big[i\in\eta^A_s]\big|
\leq\P\big[i\in\eta^{\{l\}}_s\big],
\ee
we find that
\be
2\big|\Ga(P_sf_i,P_sf_j)(A)\big|
\leq\sum_{k\in A}\sum_{l\not\in A}a(k,l)
\P\big[i\in\eta^{\{l\}}_s\big]\P\big[j\in\eta^{\{l\}}_s\big]
+\de\sum_{k\in A}\P\big[i\in\eta^{\{k\}}_s\big]\P\big[j\in\eta^{\{k\}}_s\big],
\ee
which by (\ref{varfor}) implies that
\be\ba{l}
\dis\big|\cov\big(f_i(\eta^{\{0\}}_t),f_j(\eta^{\{0\}}_t)\big)\big|\\[5pt]
\dis\quad\leq\int_0^t\sum_{k,l}a(k,l)
\P\big[k\in\eta^{\{0\}}_{t-s},\ l\not\in\eta^{\{0\}}_{t-s}\big]
\P\big[i\in\eta^{\{l\}}_s\big]\P\big[j\in\eta^{\{l\}}_s\big]\di s\\[5pt]
\dis\quad\phantom{\leq}+\de\int_0^t\sum_k\P\big[k\in\eta^{\{0\}}_{t-s}\big]
\P\big[i\in\eta^{\{k\}}_s\big]\P\big[j\in\eta^{\{k\}}_s\big]\di s.
\ec
Inserting this into (\ref{fifj}), we obtain for the quantity in (\ref{mud})
the estimate
\be\ba{l}\label{mudest}
\dis\int_0^\infty e^{-\la t}\di t\,\sum_{i,j}\psi_\ga(i,j)
\P\big[i\in\eta^{\{0\}}_t,\ j\in\eta^{\{0\}}_t\big]\\[5pt]
\dis\quad\leq\int_0^\infty e^{-\la t}\di t\,\sum_{i,j}\psi_\ga(i,j)
\P\big[i\in\eta^{\{0\}}_t\big]\P\big[j\in\eta^{\{0\}}_t\big]\\[5pt]
\dis\quad\phantom{\leq}+\int_0^\infty e^{-\la t}\di t\int_0^t\di s\,
\sum_{i,j,k,l}\psi_\ga(i,j)a(k,l)
\P\big[k\in\eta^{\{0\}}_{t-s},\ l\not\in\eta^{\{0\}}_{t-s}\big]
\P\big[i\in\eta^{\{l\}}_s\big]\P\big[j\in\eta^{\{l\}}_s\big]\\[5pt]
\dis\quad\phantom{\leq}+\de\int_0^\infty e^{-\la t}\di t\int_0^t\di s
\,\sum_{i,j,k}\psi_\ga(i,j)\P\big[k\in\eta^{\{0\}}_{t-s}\big]
\P\big[i\in\eta^{\{k\}}_s\big]\P\big[j\in\eta^{\{k\}}_s\big].
\ec
Here
\be\ba{l}
\dis\sum_{i,j,k}\psi_\ga(i,j)\P\big[k\in\eta^{\{0\}}_{t-s}\big]
\P\big[i\in\eta^{\{k\}}_s\big]\P\big[j\in\eta^{\{k\}}_s\big]\\[5pt]
\dis\quad=\sum_{i,j,k}\psi_\ga(k^{-1}i,k^{-1}j)\P\big[k\in\eta^{\{0\}}_{t-s}\big]
\P\big[k^{-1}i\in\eta^{\{0\}}_s\big]\P\big[k^{-1}j\in\eta^{\{0\}}_s\big]\\[5pt]
\dis\quad=\Big(\sum_k\P\big[k\in\eta^{\{0\}}_{t-s}\big]\Big)
\Big(\sum_{i,j}\psi_\ga(i,j)\P\big[i\in\eta^{\{0\}}_s\big]
\P\big[j\in\eta^{\{0\}}_s\big]\Big)\\[5pt]
\dis\quad=\E\big[|\eta^{\{0\}}_{t-s}|\big]
\sum_{i,j}\psi_\ga(i,j)\P\big[i\in\eta^{\{0\}}_s\big]
\P\big[j\in\eta^{\{0\}}_s\big]
\ec
and similarly
\be\ba{l}
\dis\sum_{i,j,k,l}\psi_\ga(i,j)a(k,l)
\P\big[k\in\eta^{\{0\}}_{t-s},\ l\not\in\eta^{\{0\}}_{t-s}\big]
\P\big[i\in\eta^{\{l\}}_s\big]\P\big[j\in\eta^{\{l\}}_s\big]\\[5pt]
\dis\quad\leq\sum_{i,j,k,l}\psi_\ga(l^{-1}i,l^{-1}j)a(k,l)
\P\big[k\in\eta^{\{0\}}_{t-s}\big]
\P\big[l^{-1}i\in\eta^{\{0\}}_s\big]\P\big[l^{-1}j\in\eta^{\{0\}}_s\big]\\[5pt]
\dis\quad=\Big(\sum_{k,l}a(k,l)\P\big[k\in\eta^{\{0\}}_{t-s}\big]\Big)
\Big(\sum_{i,j}\psi_\ga(i,j)
\P\big[i\in\eta^{\{0\}}_s\big]\P\big[j\in\eta^{\{0\}}_s\big]\Big)\\[5pt]
\dis\quad=|a|\,\E\big[|\eta^{\{0\}}_{t-s}|\big]\sum_{i,j}\psi_\ga(i,j)
\P\big[i\in\eta^{\{0\}}_s\big]\P\big[j\in\eta^{\{0\}}_s\big].
\ec
Inserting this into (\ref{mudest}) and recalling that this is an estimate for
the quantity in (\ref{mud}) yields
\be\ba{l}\label{dijest}
\dis\int\hat\mu_\la(\di A)1_{\{0\in A\}}e_\ga(A)\\[5pt]
\dis\quad\leq\int_0^\infty e^{-\la t}\di t\,\sum_{i,j}\psi_\ga(i,j)
\P\big[i\in\eta^{\{0\}}_t\big]\P\big[j\in\eta^{\{0\}}_t\big]\\[5pt]
\dis\quad\phantom{\leq}+(|a|+\de)\int_0^\infty e^{-\la t}\di t\int_0^t\di s\,
\E\big[|\eta^{\{0\}}_{t-s}|\big]\sum_{i,j}\psi_\ga(i,j)
\P\big[i\in\eta^{\{0\}}_s\big]\P\big[j\in\eta^{\{0\}}_s\big]\\[5pt]
\dis\quad=
\Big(1+(|a|+\de)\int_0^\infty e^{-\la t}\di t\,
\E\big[|\eta^{\{0\}}_t|\big]\Big)
\Big(\int_0^\infty e^{-\la t}\di t\,\sum_{i,j}\psi_\ga(i,j)
\P\big[i\in\eta^{\{0\}}_t\big]\P\big[j\in\eta^{\{0\}}_t\big]\Big),
\ec
where in the last step we have changed the integration order on the set
$\{(s,t):0\leq s\leq t\}$. Using the fact that $\psi_\ga(i,j)=e^{\ga d(i,j)}$
where $d$ is a metric, we may further estimate the sum in the second factor on
the right-hand side of (\ref{dijest}) as
\be\ba{l}
\dis\sum_{i,j}\psi_\ga(i,j)
\P\big[i\in\eta^{\{0\}}_t\big]\P\big[j\in\eta^{\{0\}}_t\big]
=\sum_{i,j}e^{\ga d(i,j)}
\P\big[i\in\eta^{\{0\}}_t\big]\P\big[j\in\eta^{\{0\}}_t\big]\\[5pt]
\leq\dis\quad\sum_{i,j}e^{\ga(d(0,i)+d(0,j))}
\P\big[i\in\eta^{\{0\}}_t\big]\P\big[j\in\eta^{\{0\}}_t\big]\\[5pt]
=\dis\quad\big(\sum_ie^{\ga d(0,i)}\P\big[i\in\eta^{\{0\}}_t\big]\big)^2
=\E\big[\sum_{i\in\eta^{\{0\}}_t}e^{\ga d(0,i)}\big]^2.
\ec
Inserting this into (\ref{dijest}) and recalling the definition of
$\hat\pi_\la$ in (\ref{hatpi}) yields
\be\label{appi}
\int\hat\mu_\la(\di A)1_{\{0\in A\}}e_\ga(A)
\leq \Big(1+(|a|+\de)\hat\pi_\la\Big)\int_0^\infty e^{-\la t}\di t\,
\E\big[e_\ga(\eta^{\{0\}}_t)\big]^2.
\ee
We note that setting $\ga=0$ in (\ref{rga}) shows that
\be\label{erE}
e^{rt}\leq\E\big[|\eta^{\{0\}}_t|\big]\qquad(t\geq 0),
\ee
and therefore
\be
\lim_{\la\down r}\hat\pi_\la=
\lim_{\la\down r}\int_0^\infty e^{-\la t}\di t\,\E\big[|\eta^{\{0\}}_t|\big]
=\infty.
\ee
Using this and (\ref{appi}), we arrive at (\ref{unexmom}).\qed

\noi
As a direct applications  we obtain:\med

\noi
\tbf{Proof of Lemma~\ref{L:unmom}} This is special case of
Lemma~\ref{L:unexmom}, where $\ga=0$.\qed

\noi
\tbf{Proof of Lemma~\ref{L:nutight}}
This is very similar to the proof of Lemma~\ref{L:finex}.
For $\de\in(\de_{\rm c},\infty)$, let $(\eta^{\de,\{0\}}_t)_{t\geq 0}$ and
$\nucirc_\de$ be as in Lemma~\ref{L:nutight}. Let $\La_k$ be finite sets such
that $0\in\La_k\sub\La$ and $\La_k\up\La$. It is again easy to check that
$A\mapsto f_k^{\gamma}(A):=e_\ga(A\cap\La_k)1_{\{0\in A\}}$ is a continuous,
compactly supported real function on $\Pc_+$. Therefore, since (by
Proposition~\ref{P:conv}) the $\frac{1}{\hat\pi_{\la_n}}\hat\mu_{\la_n}$
converge vaguely to $\nucirc^{\de}$,
\begin{eqnarray*}
\int\nucirc^{\de}(\di A)f_k^{\ga}(A)
=\lim_{n\to\infty}
\frac{1}{\hat\pi_{\la_n}}\int\hat\mu_{\la_n}(\di A)f_k^{\ga}(A)
&\leq&\liminf_{n\to\infty}
\frac{1}{\hat\pi_{\la_n}}\int\hat\mu_{\la_n}(\di A) e_{\ga}(A)1_{\{0\in A\}}\\
&\leq&(|a|+\de)\int_0^\infty e^{-rt}\di t\,
\E\big[e_\ga(\eta^{\de,\{0\}}_t)\big]^2.
\end{eqnarray*}
Letting $k\up\infty$ such that $f_k^{\gamma}\uparrow e_{\ga}(A)1_{\{0\in A\}}$
we arrive at (\ref{nutight}) by the monotone convergence theorem.
\qed

\appendix

\section{Exponential decay in the subcritical regime}\label{A:decay}

\subsection{Statement of the result}

The aim of this appendix is to show how the arguments in \cite{AJ07}, which
are written down for contact processes on transitive graphs, can be extended
to prove Theorem~\ref{T:rprop}~(d) for the class of $(\La,a,\de)$-contact
processes considered in this article. To formulate this properly, only in this
appendix, we will consider a class of contact processes that is more general
than both the one defined in Section~\ref{S:group} and the one considered in
\cite{AJ07}, and contains them both as subclasses. Indeed, only in this
appendix, will we drop the assumptions that $\La$ has a group structure (as in
the rest of this article) or that $\La$ has a graph structure (as in
\cite{AJ07}). The only structure on $\La$ that we will use is the structure
given by the infection rates $(a(i,j))_{i,j\in\La}$.

Let $\La$ be any countable set and let $a:\La\times\La\to\half$ be a
function. By definition, an \emph{automorphism} of $(\La,a)$ is a bijection
$g:\La\to\La$ such that $a(gi,gj)=a(i,j)$ for each $i,j\in\La$. Let ${\rm
Aut}(\La,a)$ denote the group of automorphisms of $(\La,a)$. We say that a
subgroup $G\sub{\rm Aut}(\La,a)$ is \emph{(vertex) transitive} if for each
$i,j\in\La$ there exists a $g\in G$ such that $gi=j$. In particular, we say
that $(\La,a)$ is transitive if ${\rm Aut}(\La,a)$ is transitive.

Let $(\La,a)$ be transitive, let $a^\dgg(i,j):=a(j,i)$, and assume that
\be\label{genassum}
|a|:=\sum_{j\in\La}a(i,j)<\infty
\quad\mbox{and}\quad
|a^\dgg|:=\sum_{j\in\La}a^\dgg(i,j)<\infty,
\ee
where by the transitivity of $(\La,a)$, these definitions do not depend on the
choice of $i\in\La$. Then, for each $\de\geq 0$, there exists a well-defined
contact process on $\La$ with generator as in (\ref{Gdef}) and also the dual
contact process with $a$ replaced by $a^\dgg$ is well-defined. \emph{Only in
this appendix}, we will use the term $(\La,a,\de)$-contact process
(resp.\ $(\La,a^\dgg,\de)$-contact process) in this more general sense.

For any $(\La,a,\de)$-contact process, as defined in this appendix, we define
the critical recovery rate $\de_{\rm c}=\de_{\rm c}(\La,a)$ as in (\ref{dec}),
which satisfies $\de_{\rm c}<\infty$ but may be zero in the generality
considered here. A straightforward extension of \cite[Lemma~1.1]{Swa09} shows
that the exponential growth rate $r=r(\La,a,\de)$ in (\ref{rdef}) is
well-defined for the class of $(\La,a,\de)$-contact processes considered
here.

We will show that the arguments in \cite{AJ07} imply the following result.
\begin{atheorem}
\pbf{(Exponential decay in the subcritical regime)}\label{T:subex}
Let $(\La,a)$ be transitive and let $a$ satisfy (\ref{genassum}). Then
$\{\de\geq 0:r(\La,a,\de)<0\}=(\de_{\rm c},\infty)$.
\end{atheorem}
We remark that Theorem~\ref{T:rprop}~(a) does not hold in general for the
class of $(\La,a,\de)$-contact processes considered in this appendix. This is
related to unimodularity. A transitive subgroup $G\sub{\rm Aut}(\La,a)$ is
\emph{unimodular} if \cite[formula~(3.3)]{BLPS99}
\be\label{unimod}
|\{gi:g\in G,\ gj=j\}|=|\{gj:g\in G,\ gi=i\}|\qquad(i,j\in\La).
\ee
Note that this is trivially satisfied if $\La$ is a group and $G=\La$ acts on
itself by left multiplication, in which case the sets on both sides of the
equation consist of a single element. Unimodularity gives rise to the {\em
mass transport principle} which says that for any function
$f:\La\times\La\to\half$ such that $f(gi,gj)=f(i,j)$ $(g\in G,\ i,j\in\La)$,
one has $\sum_jf(i,j)=\sum_jf(j,i)$. In particular, this implies that the
constants $|a|$ and $|a^\dgg|$ from (\ref{genassum}) are equal and that
$r(\La,a,\de)=r(\La,a^\dgg,\de)$. In the nonunimodular case, this is in
general no longer true and in fact it is not hard to construct examples where
the critical recovery rates $\de_{\rm c}(\La,a)$ and $\de_{\rm c}(\La,a^\dgg)$
of a contact process and its dual are different. We remark that although in
\cite{AJ07}, the authors do not always clearly distinguish between a contact
process and its dual (e.g., in their formulas (1.3), (1.9) and Lemma~1.4),
they do not assume that $a=a^\dgg$ and their results are valid also in the
asymmetric case $a\neq a^\dgg$.

\subsection{The key differential inequalities and their consequences}

The main method used in \cite{AJ07}, that in its essence goes back to
\cite{AB87} and that yields Theorem~\ref{T:subex} and a number of related
results, is the derivation of differential inequalities for certain quantities
related to the process. Using the graphical representation to construct a
$(\La,a,\de)$-contact process and its dual, we define the \emph{susceptibility}
as
\be
\chi=\chi(\La,a,\de)=\E\big[\int_0^\infty\!|\eta^{\{0\}}_t|\,\di t\big],
\ee
which may be $+\infty$. Moreover, letting $\poi^{\rm c}$ be a Poisson point
process on $\La\times\R$ with intensity $h\geq 0$, independent of the Poisson
point processes $\poi^{\rm i}$ and $\poi^{\rm r}$ corresponding to infection
arrows and recovery symbols, we define
\be
\tet=\tet(\La,a,\de,h):=\P\big[C_{(0,0)}\cap\poi^{\rm c}\neq\emptyset\big]
\quad\mbox{where}\quad
C_{(i,s)}:=\big\{(j,t):t\geq s,\ (i,s)\leadsto(j,t)\big\}.
\ee
Then $\tet$ can be interpreted as the density of infected sites in the upper
invariant law of a (dual) ``$(\La,a^\dgg,\de,h)$-contact process'', which in
addition to the dynamics in (\ref{Gdef}) exhibits spontaneous infection of
healthy sites with rate $h$, corresponding to a term in the generator of the
form $h\sum_i\{f(A\cup\{i\})-f(A)\}$.

Let $\La,a,\de$ be fixed and for $\la,h\geq 0$ let
$\tet=\tet(\la,h):=\tet(\La,\la a,\de,h)$ and $\chi=\chi(\la):=\chi(\La,\la
a,\de)$ be the quantities defined above. The analysis in \cite{AJ07} centers
on the deriviation of the following three differential inequalities (see
\cite[formulas (1.17), (1.19) and (1.20)]{AJ07})
\be\ba{rr@{\,}c@{\,}l}\label{keydif}
{\rm(i)}&\dis\dif{\la}\chi&\leq&\dis|a|\chi^2,\\[5pt]
{\rm(ii)}&\dis\dif{\la}\tet&\leq&\dis|a|\tet\dif{h}\tet,\\[5pt]
{\rm(iii)}&\dis\tet&\leq&\dis h\dif{h}\tet
+\big(2\la^2|a|\tet+h\la\big)\dif{\la}\tet+\tet^2.
\ec
These differential inequalities, and their proofs, generalize without a change
to the more general class of $(\La,a,\de)$-contact processes discussed in this
appendix.

Since $\tet\geq h(1+h)$, which follows by estimating the $(\La,\la
a^\dgg,\de,h)$-contact process from below by a process with no infections, one
has $h\leq\tet(1-\tet)$. Inserting this into (\ref{keydif})~(iii) yields
\be\label{key2}
\tet\leq h\dif{h}\tet+\Big(2\la^2|a|
+\frac{\la}{1-\tet}\Big)\tet\dif{\la}\tet+\tet^2.
\ee
Abstract results of Aizenman and Barsky \cite[Lemmas~4.1 and 5.1]{AB87} allow
one to draw the following conclusions from (\ref{keydif})~(ii) and
(\ref{key2}).

\bl\pbf{(Estimates on critical exponents)}\label{L:critexp}
Assume that there exists some $\la'>0$ such that $\tet(\la',0)=0$ and
$\lim_{h\to 0}h^{-1}\tet(\la',h)=\infty$. Then there exist $c_1,c_2>0$ such
that
\be\ba{rr@{\,}c@{\,}ll}\label{critexp}
{\rm(i)}&\dis\tet(\la',h)&\geq&\dis c_1h^{1/2}\qquad&(h\geq 0),\\[5pt]
{\rm(ii)}&\dis\tet(\la,0)&\geq&\dis c_2(\la-\la')\qquad&(\la\geq\la').
\ec
\el
Note that this lemma (in particular, formula (\ref{critexp})~(i), which
depends on the assumption that $\lim_{h\to 0}h^{-1}\tet(\la',h)=\infty$)
implies in particular that if for some fixed $\la'>0$, one has
$\tet(\la',h)\sim h^\al$ as $h\to 0$, then either $\al\leq\frac{1}{2}$ or
$\al\geq 1$.\med

\noi
\tbf{Remark} Lemmas~4.1 and 5.1 of \cite{AB87} are also cited in
\cite[Thm.~4.1]{AJ07}, but there the statement that $c_1,c_2>0$ is erroneously
replaced by the (empty) statement that $c_1,c_2<\infty$.\med

\noi
\tbf{Proof of Theorem~\ref{T:subex} (sketch)} Set
\bc
\dis\la_{\rm c}&:=&\dis\inf\{\la\geq 0:\tet(\la,0)>0\},\\[5pt]
\dis\la'_{\rm c}&:=&\dis\inf\{\la\geq 0:\chi(\la)=\infty\}.
\ec
Since $\chi(\la)<\infty$ implies $\tet(\la,0)=0$, obviously $\la'_{\rm
c}\leq\la_{\rm c}$. Our first aim is to show that they are in fact equal. We
note that it is always true that $\la'_{\rm c}>0$. It may happen that
$\la'_{\rm c}=\infty$ but in this case also $\la_{\rm c}=\infty$ so without
loss of generality we may assume that $\la'_{\rm c}<\infty$.

It follows from (\ref{keydif})~(i) and approximation of infinite systems by
finite systems  (compare \cite[Lemma~3.1]{AN84}, which is written down for
unoriented percolation and which is cited in \cite[formula~(1.18)]{AJ07}) that
$\lim_{\la\up\la'_{\rm c}}\chi(\la)=\chi(\la'_{\rm c})=\infty$, and in fact
\be\label{chicrit}
\chi(\la)\geq\frac{|a|^{-1}}{\la'_{\rm c}-\la}\qquad(\la<\la'_{\rm c}).
\ee
Now either $\tet(\la'_{\rm c},0)>0$, in which case we are done, or
$\tet(\la'_{\rm c},0)=0$. In the latter case, since
\be\label{diftet}
\chi(\la)=\lim_{h\to 0}h^{-1}\tet(\la,h)\qquad(\la<\la'_{\rm c}),
\ee
(see \cite[formula~(1.11)]{AJ07}), using the monotonicity of $\tet$ in $\la$
and $h$, it follows from (\ref{chicrit}) that
\be
\lim_{h\to 0}h^{-1}\tet(\la'_{\rm c},h)=\infty
\ee
and therefore Lemma~\ref{L:critexp} implies that (\ref{critexp}) holds at
$\la'=\la'_{\rm c}$. In particular, (\ref{critexp})~(ii) implies that
$\tet(\la,0)>0$ for $\la>\la'_{\rm c}$, hence $\la_{\rm c}=\la'_{\rm c}$.

Since by a trivial rescaling of time, questions about critical values for
$\la$ can always be translated into questions about critical values for $\de$,
we learn from this that for any $(\La,a,\de)$-contact process, one has
$\chi(\La,a,\de)<\infty$ if $\de>\de_{\rm c}(\La,a)$, where the latter
critical point is defined in (\ref{dec}). It follows from (\ref{erE}) that
$\chi(\La,a,\de)=\infty$ if $r(\de)=r(\La,a,\de)\geq 0$, hence we must have
$r(\de)<0$ for $\de\in(\de_{\rm c},\infty)$. Part~(b) of Theorem~\ref{T:rprop}
is easily generalized to the class of $(\La,a,\de)$-contact processes
considered in this appendix. Moreover, it is not hard to prove that $r<0$
implies that the process does not survive. This shows that $r(\de)\geq 0$ on
$[0,\de_{\rm c})$ while $\de\mapsto r(\de)$ is continuous, which allows us to
conclude that $\{\de\geq 0:r(\de)<0\}=(\de_{\rm c},\infty)$ if $\de_{\rm
c}>0$. If $\de_{\rm c}=0$ (which may happen for the general class of models
considered here), then we may use the fact that $\tet(\La,a,0)=1$ to conclude
that $r(\La,a,0)\geq 0$, hence the conclusion of Theorem~\ref{T:subex} is also
valid in this case.\qed

\section{Some results on quasi-invariant laws}\label{A:quasi}

In this appendix we collect some basic results on $\la$-positivity and
quasi-invariant laws for which we did not find an exact reference in the
literature. We will be interested in continuous-time Markov chains taking
values in a countable set $S$, which may have a finite lifetime due to killing
or explosion. To formalize this, let $\ov S:=S\cup\{\infty\}$ be the set $S$
with one extra point added and let $\ov Q$ be a \emph{Q-matrix} on $\ov S$,
i.e., $\ov Q:\ov S^2\to\R$ is a function such that
\be
\ov Q(i,j)\geq 0\quad(i\neq j),
\qquad\sum_{k\in\ov S,\ k\neq i}\ov Q(i,k)<\infty,
\qquad\sum_{k\in\ov S}\ov Q(i,k)=0
\ee
for all $i,j\in\ov S$. We set $\ov Q(i):=-\ov Q(i,i)$ and call $T:=\{i\in\ov
S:\ov Q(i)=0\}$ the set of \emph{traps}. We assume that $\infty$ is a trap,
i.e, $\ov Q(\infty)=0$. For any initial law on $\ov S$, we construct a
continuous-time Markov chain $X=(X_t)_{t\geq 0}$ with Q-matrix $\ov Q$ in the
usual way from its embedded Markov chain. More precisely, let $Y=(Y_k)_{0\leq
  k<N+1}$ be a Markov chain in $\ov S$ with possibly finite lifetime
$N=\inf\{k\geq 0:Y_k\in T\}$ and transition probabilities
\be
\P[Y_k=j\,|\,Y_{k-1}=i]=\ov Q(i,j)/\ov Q(i)\qquad(0<k<N+1).
\ee
Conditional on $Y$, let $\sig_k$ be independent, exponentially
distributed random variables with parameter $\ov Q(Y_k)$ $(0\leq k<N+1)$
(in particular, $\sig_N=\infty$ if $N<\infty$), set
$\tau_n:=\sum_{0\leq k<n}\sig_k$ $(0\leq n\leq N+1)$, and let $\tau:=\tau_{N+1}$,
which may be finite if $N=\infty$. Then setting
\be
X_t:=Y_k\quad(\tau_k\leq t<\tau_{k+1})
\ee
defines a continuous-time Markov chain $X=(X_t)_{0\leq t<\tau}$ with Q-matrix
$\ov Q$ and possible finite lifetime $\tau$. We call $\{\tau<\infty\}$ the
event of \emph{explosion} and say that the process is \emph{nonexplosive} if
this has probability zero. By construction, each trap $i\in T$ has the
property that $X_s=i$ for some $s\in[0,\tau)$ implies $X_t=i$ for all
$t\in[s,\tau)$. In particular, this is true for $i=\infty$ which is a trap by
assumption.

We will only be interested in the process $X$ as long as it stays in $S$. If
$X$ jumps to $\infty$ at some point, then we say that the process gets
\emph{killed}. We call $\ov Q(i,\infty)$ the \emph{killing rate} at~$i$. If
$\ov Q(i,\infty)=0$ for all $i\in S$ then we say the process has \emph{zero
  killing rates}. If the process explodes, then we also set $X_t:=\infty$ for
all $t\geq\tau$, i.e., we use the same cemetery state $\infty$
regardless of whether the process disappears from $S$ due to it being killed
or due to explosion.

We let
\be
\ov P_t(i,j):=\P[X_{s+t}=j\,|\,X_s=i]\quad(s,t\geq 0,\ i,j\in\ov S)
\ee
denote the transition probabilities of the Markov process $X$ in $\ov S$,
and let $P_t$ denote the restriction of $\ov P_t$ to $S^2$. 
Due to the possibility of killing or explosion, the $(P_t)_{t\geq 0}$ are in
general subprobability kernels on $S$. Let $Q$ denote the restriction
of $\ov Q$ to $S^2$. It follows from well-known results (see
e.g.\ \cite[Prop~2.30]{Lig10}, \cite[Thm~2.8.4]{Nor97}) that the functions
$t\mapsto P_t(i,j)$ are continuously differentiable for each $i,j\in S$ and
that the $(P_t)_{t\geq 0}$ are given by the minimal nonnegative solution to
the Kolmogorov backward equations
\be\label{difp}
\dif{t}P_t(i,k)=\sum_{j\in S}Q(i,j)P_t(j,k)\quad(t\geq 0,\ i,k\in S)
\ee
with initial condition $P_0(i,j)=1_{\{i=j\}}$. We say that the process is {\em
 irreducible on $S$} if for each $i,j\in S$ there exist $i=i_0,\ldots,i_n=j$
such that $Q(i_{k-1},i_k)>0$ for $k=1,\ldots,n$. If the process is irreducible
on $S$, then by \cite[Thm~1]{Kin63} the \emph{decay parameter}
\be\label{laS}
\la_S:=-\lim_{t\to\infty}t^{-1}\log P_t(i,j)
\ee
exists and does not depend on $i,j\in S$. For nonexplosive processes with zero
kiling rates, we define transience, null-recurrence and positive recurrence in
the standard way. We use the usual notation for matrices and vectors
indexed by $S$, i.e., $AB(i,k):=\sum_jA(i,j)B(j,k)$, $gA(i):=\sum_jg(j)A(j,i)$
and $Ah(i):=\sum_jA(i,j)h(j)$. With these definitions, one has the following
facts that will be proven below.

\bl\pbf{(Doob transformed process)}\label{L:htrafo}
Assume that $h:S\to(0,\infty)$ and $\la\in\R$ satisfy $P_th=e^{-\la t}h$
$(t\geq 0)$. Then
\be\label{Phdef}
P^h_t(i,j):=e^{\la t}h(i)^{-1}P_t(i,j)h(j)\qquad(i,j\in S,\ t\geq 0)
\ee
are the transition probabilities of a nonexplosive continuous-time Markov
chain in $S$ with zero killing rates and with Q-matrix given by
\be\label{Qhdef}
Q^h(i,j):=h(i)^{-1}Q(i,j)h(j)+\la 1_{\{i=j\}}\qquad(i,j\in S).
\ee
In particular, $\sum_{j:\,j\neq i}Q^h(i,j)<\infty$ and $\sum_jQ^h(i,j)=0$ for
all $i\in S$.
\el

\bl\pbf{($\la$-positivity)}\label{L:Rpos}
Assume that $Q$ is irreducible on $S$ and that $g,h:S\to(0,\infty)$ and
$\la\in\R$ satisfy
\be\label{ghcond}
gP_t=e^{-\la t}g,\quad P_th=e^{-\la t}h\quad(t\geq 0)
\quand c:=\sum_ig(i)h(i)<\infty.
\ee
Then the transition probabilities $(P^h_t)_{t\geq 0}$ in (\ref{Phdef}) belong
to a positively recurrent continuous-time Markov chain with unique invariant
law given by $\pi(i)=c^{-1}g(i)h(i)$. Moreover, the conditions (\ref{ghcond})
determine $g$ and $h$ uniquely up to multiplicative constants and imply that
$\la=\la_S$.
\el

\bl\pbf{(Quasi-invariant law)}\label{L:quasi}
In the set-up of Lemma~\ref{L:Rpos}), assume moreover that $\inf_{i\in
 S}h(i)>0$. Then the process $X$ started in any deterministic initial state
$i\in S$ satisfies
\be\label{ratlim}
\P^i[X_t\in\,\cdot\,|X_t\neq\infty]\Asto{t}\nu,
\ee
where $\nu$ is the probability measure on $S$ defined by $\nu(i):=g(i)/\sum_j
g(j)$ $(i\in S)$ and $\Rightarrow$ denotes weak convergence of probability
measures on $S$.
\el

Before we sketch the proofs of these results, we first discuss what can be
found about this in the literature. For discrete-time Markov chains and
more generally for countable nonnegative matrices, the concepts of
R-transience, R-null recurrence, and R-positivity were introduced by
Vere-Jones \cite{Ver67} (which builds on his D.Phil.\ thesis from
1961). Kingman \cite{Kin63} then treated the continuous-time case. (A good
general reference to this material is \cite[Sect~5.2]{And91}.) In a more
general set-up than ours, Kingman proved that the limit in (\ref{laS}) exists.
He then defined $(P_t)_{t\geq 0}$ to be $\la$-transient or $\la$-recurrent
depending on whether
\be
\int_0^\infty P_t(i,i)e^{\la_St} dt<\infty\quad\mbox{or}\quad=\infty,
\ee
where as a result of irreducibility the definition does not depend on the
choice of the reference point $i\in S$. In the $\la$-recurrent case, he called
$(P_t)_{t\geq 0}$ $\la$-null-recurrent or $\la$-positive (recurrent) depending
on whether
\be
\lim_{t\to\infty}P_t(i,i)e^{\la_St}=0\quad\mbox{or}\quad>0,
\ee
where the limit is shown to exist and the definition does not depend on the
choice of $i\in S$. In the $\la$-recurrent case, he showed \cite[Thm~4]{Kin63}
that there are functions $g,h:S\to(0,\infty)$, unique up to multiplicative
constants, such that $gP_t=e^{-\la_S t}g$, $P_th=e^{-\la_S t}h$ $(t\geq 0)$,
and these satisfy $\sum_ig(i)h(i)<\infty$ if and only if $(P_t)_{t\geq 0}$ is
$\la$-positive. He moreover defined $(P^h_t)_{t\geq 0}$ as in (\ref{Phdef})
and observed that these are the transition probabilities of a continuous-time
Markov process.\footnote{In fact, Kingman defines a (right) $\la$-subinvariant
  vector to be any function $h:S\to(0,\infty)$ such that $P_th\leq e^{-\la_S
    t}h$ $(t\geq 0)$. He proves that such $\la$-subinvariant vectors exist
  quite generally and defines $(P^h_t)_{t\geq 0}$ as in (\ref{Phdef}) for any
  such $h$, which may now be subprobability kernels corresponding to a process
  with killing.}

Some care is needed in applying Kingman's results to our setting, however, since
his setting is more general than ours. He assumes that $S$ is an irreducible
subclass of some larger space $\ov S$ that may be more complicated than in our
setting, and he only assumes that the Markov process corresponding to
$(\ov P_t)_{t\geq 0}$ a.s.\ assumes values in the countable set $\ov S$ at
deterministic times. This includes processes that are not defined by a
Q-matrix and that leave each state instantaneously, such as Blackwell's
example \cite[Sect.~2.4]{Lig10} or the FIN diffusion defined in \cite{FIN02},
which is a Brownian motion time-changed in such a way that at deterministic
times it a.s.\ takes values in a countable, dense subset of the real line.
While Kingman's set-up is in all respects more general than ours, his results
are also weaker in some respects since he does not prove that the transformed
transition probabilities $(P^h_t)_{t\geq 0}$ as in (\ref{Phdef}) come from a
Q-matrix.

In practical situations, one often does not have a direct way of verifying
that a function $h$ satisfies $P_th=e^{-\la t}h$ $(t\geq 0)$, but instead
starts off from a solution to $Qh=-\la h$. The latter equation is in general
not enough to guarantee the first one, so extra conditions are needed, see
\cite{NP93}. In our case, however, solutions to $P_th=e^{-\la t}h$ can be
obtained directly from the eigenmeasures, which is why the lemmas above are
sufficient for our purposes.

Probability measures $\nu$ satisfying $\nu P_t=e^{-\la_St}\nu$ are called
\emph{quasi-invariant laws} and (\ref{ratlim}) is a \emph{ratio limit
 theorem}. We refer to \cite{FKM96} and references therein for a more
detailed discussion of these concepts.

To prepare for the proofs of Lemmas~\ref{L:htrafo}--\ref{L:quasi}, we prove one
technical lemma.

\bl\pbf{(Cadlag processes have well-defined rates)}\label{L:cadlagQ}
Let $S$ be a countable set and let $(P_t)_{t\geq 0}$ be probability kernels on
$S$ such that $P_sP_t=P_{s+t}$ $(s,t\geq 0)$ and $\lim_{t\down
 0}P_t(i,i)=P_0(i,i)=1$ $(i\in S)$. Assume that for each $i\in S$, there
exists a Markov process $X=(X_t)_{t\geq 0}$ in $S$ with initial state $X_0=i$,
transition probabilities $(P_t)_{t\geq 0}$, and cadlag sample paths. Then
there exists a Q-matrix on $S$ such that $(P_t)_{t\geq 0}$ is the minimal
nonnegative solution of (\ref{difp}).
\el
\tbf{Proof} Define inductively stopping times by
$\tau_0=\tau^\eps_0=0$ and
\bc
\dis\tau_k&:=&\dis\inf\{t\geq\tau_{k-1}:X_t\neq X_{\tau_{k-1}}\}\\[5pt]
\dis\tau^\eps_k&:=&\dis\inf\{\eps l\geq\tau^\eps_{k-1}:
X_{\eps l}\neq X_{\tau^\eps_{k-1}},\ l\in\N\}\qquad(\eps>0).
\ec
Let $N_\eps:=1+\sup\{k\geq 0:\tau^\eps_k<\infty\}$. Then, for each $\eps>0$,
we may define a Markov chain $Y^\eps=(Y^\eps_k)_{0\leq k<N_\eps+1}$ by
$Y^\eps_k:=X_{\tau^\eps_k}$ $(0\leq k<N_\eps+1)$. Conditional on $Y^\eps$,
the holding times $(\tau^\eps_{k+1}-\tau^\eps_k)$ with $(0\leq k<N_\eps+1)$
are independent and geometrically distributed. By the fact that $X$ has cadlag
sample paths, $Y^\eps\to Y$ a.s.\ where the embedded Markov chain
$Y=(Y_k)_{0\leq k<N+1}$ is defined analogously to $Y^\eps$ with $\tau^\eps_k$
replaced by $\tau_k$. Moreover, the collection of times $(\tau^\eps_k)_{0\leq
  k<N_\eps+1}$ a.s.\ converges to $(\tau_k)_{0\leq k<N+1}$.

In particular, for the process started in $i$, since $\tau^\eps_1$ is
geometrically distributed and $\tau^\eps_1\to\tau_1$ a.s., we see that
$\tau_1$ is exponentially distributed and the limit
\be
Q(i)=-Q(i,i):=\lim_{\eps\down 0}\eps^{-1}\big(1-P_\eps(i,i)\big)
\ee
exists, where we allow for the case $Q(i)=0$ (which corresponds to
$\tau_1=\infty$). (Note that the assumption of cadlag sample paths implies
$\tau_1>0$.) If $Q(i)>0$, then we observe that $Y^\eps_1$ is distributed
according to the law
\be
\P^i\big[Y^\eps_1=j\big]=\big(1-P_\eps(i,i)\big)^{-1}P_\eps(i,j)
\qquad(j\in S,\ j\neq i).
\ee
Since $Y^\eps_1\to Y_1$ as $\eps\to 0$, we conclude that the limit
\be
Q(i,j):=\lim_{\eps\down 0}\eps^{-1}P_\eps(i,j)=Q(i)\P^i[Y_1=j]\qquad(i\neq j)
\ee
exists and satisfies $\sum_{j:\,j\neq i}Q(i,j)<\infty$ and $\sum_jQ(i,j)=0$.

It is now not hard to check that $Y$ is a Markov chain that jumps from a state
$i$ with $Q(i)>0$ to a state $j$ with probability $Q(i)^{-1}Q(i,j)$, and that
conditional on $Y$, the times $(\tau_{k+1}-\tau_k)$ are independent and
exponentially distributed with parameter $Q(Y_k)$. By
\cite[Thm~2.8.4]{Nor97}, we conclude that $(P_t)_{t\geq 0}$ is the unique
minimal nonnegative solution of (\ref{difp}).\qed

\noi
\tbf{Proof of Lemma~\ref{L:htrafo}} The fact that $P_th=e^{-\la t}h$ implies
that the $(P^h_t)_{t\geq 0}$ are probability kernels satisfying
$P^h_sP^h_t=P^h_{s+t}$ $(s,t\geq 0)$ and $\lim_{t\down
 0}P^h_t(i,i)=P^h_0(i,i)=1$ $(i\in S)$. It is not immediately clear, however,
that these are the transition probabilities of a Markov process with Q-matrix
as in (\ref{Qhdef}), or, in fact, that the latter is even well-defined.

To prove this, view $\ov S$ as the one-point compactification of $S$ (with the
appropriate topology). Then the process $X$ has sample paths in the space of
those cadlag functions $\om:\half\to\ov S$ for which $\om_s=\infty$ implies
$\om_s=\infty$ for all $t\geq s$. We may construct $X$ in the canonical way on
this space, where $X_t(\om)=\om_t$ is the coordinate projection. Let $\P^i$ be
the law of the process started in $i\in S$. We may consistently define a
new probability law $\P^{h,i}$ by
\be
\P^{h,i}\big[(X_s)_{0\leq s\leq t}\in\di\om\big]
:=e^{\la t}1_{\{\om_s\in S\ \forall 0\leq s\leq t\}}
\frac{h(\om_t)}{h(i)}\P^i\big[(X_s)_{0\leq s\leq t}\in\di\om\big]
\qquad(t\geq 0).
\ee
Then $P^{h,i}$ is concentrated on cadlag paths $\om:\half\to S$.
It is straightforward to check that under this new law, $X$ is a Markov
process with transition kernels $(P^h_t)_{t\geq 0}$. Since $X$ has cadlag
sample paths, we may invoke Lemma~\ref{L:cadlagQ} to conclude that
$X$ is a continuous-time Markov chain with Q-matrix given by
the right-hand derivative
\be
\dif{t}P^h_t(i,j)\big|_{t=0}
=\dif{t}\big(e^{\la t}h(i)^{-1}P_t(i,j)h(j)\big)\big|_{t=0}
=h(i)^{-1}Q(i,j)h(j)+\la 1_{\{i=j\}}=Q^h(i,j).
\ee
\qed

\noi
\tbf{Proof of Lemma~\ref{L:Rpos}} Since $c:=\sum_ig(i)h(i)<\infty$, we may
define a probability law $\pi$ on $S$ by $\pi(i):=c^{-1}g(i)h(i)$. Then
\be
\pi P^h_t(i)=\sum_j c^{-1}g(j)h(j)e^{\la t}h(j)^{-1}P_t(j,i)h(i)
=c^{-1}g(i)h(i)=\pi(i)\qquad(i\in S,\ t\geq 0),
\ee
where we have used that $gP_t=e^{-\la t}g$. It follows that $\pi$ is an
invariant law for the irreducible continuous-time Markov chain with Q-matrix
as in (\ref{Qhdef}), and hence the latter is positively recurrent.
In particular,
\be
\lim_{t\to\infty}e^{\la t}P_t(i,i)=\lim_{t\to\infty}P^h_t(i,i)=\pi(i)>0
\ee
which shows that $\la=\la_S$ and thus also that $(P_t)_{t\geq 0}$ is $\la$-positive. Hence, by 
applying \cite[Thm~4]{Kin63} we obtain that $g,h:S\to(0,\infty)$ are unique
up to multiplicative constants.\qed

\noi
\tbf{Proof of Lemma~\ref{L:quasi}} Since $\inf_{i\in S}h(i)>0$, (\ref{ghcond})
implies $\sum_jg(j)<\infty$, so $\nu$ is well-defined. Moreover, for any
bounded function $f:S\to\R$, 
\be\ba{l}
\dis\E^i\big[f(X_t)\,\big|\,X_t\neq\infty\big]
=\frac{\sum_jP_t(i,j)f(j)}{\sum_jP_t(i,j)}
=\frac{e^{-\la t}h(i)\sum_jP^h_t(i,j)h(j)^{-1}f(j)}
{e^{-\la t}h(i)\sum_jP^h_t(i,j)h(j)^{-1}}\\[5pt]
\dis\quad\asto{t}
\frac{\sum_j\pi(j)h(j)^{-1}f(j)}{\sum_j\pi(j)h(j)^{-1}}
=\sum_j\nu(j)f(j),
\ec
where all sums run over $j\in S$ and we have used the ergodicity of the
positively recurrent Markov process with transition probabilities
$(P^h_t)_{t\geq 0}$ and invariant law $\pi(i)=c^{-1}g(i)h(i)$, as well as the
fact that $h^{-1}f$ and $h^{-1}$ are bounded functions by our assumption that
$\inf_{i\in S}h(i)>0$.\qed

\noi \emph{Acknowledgement} We thank the referee who handled the first versions
of this paper for two impressive referee reports, which not only found a
mistake in the original proofs but mainly greatly helped improve the
presentation. We also thank Phil Pollett for helping us find our way in the
literature concering quasi-invariant laws.

\vfill

\hspace{1cm}
\parbox[t]{6cm}{\small
Anja Sturm\\
Goldschmidtstrasse 7\\
37077 G\"ottingen\\
Germany\\
email: asturm@math.uni-goettingen.de}
\hspace{1.8cm}
\parbox[t]{6cm}{\small
Jan M.~Swart\\
Institute of Information Theory and Automation of
the ASCR (\' UTIA)\\
Pod vod\'arenskou v\v e\v z\' i 4\\
18208 Praha 8\\
Czech Republic\\
e-mail: swart@utia.cas.cz}

\end{document}